\documentclass[a4j,12pt]{article} 
\usepackage{amsmath,amssymb}
\newcommand{\qed}{\hfill \hbox{\rule[-2pt]{4pt}{7pt}}}

\newtheorem{lem}{Lemma}[section]
\newtheorem{rem}[lem]{Remark}
\newtheorem{thm}[lem]{Theorem}  

\newtheorem{cor}[lem]{Corollary}
\newtheorem{propos}[lem]{Proposition}   

\newcommand{\trace}{{\rm Tr}}    
\newcommand{\tr}{{\rm tr}}    
\newcommand{\B}{{\Bbb B}}
\newcommand{\G}{{\Bbb G}}
\newcommand{\Q}{{\Bbb Q}}
\newcommand{\Z}{{\Bbb Z}}

\newcommand{\C}{{\Bbb C}}

\newcommand{\bH}{{\Bbb H}}

\newcommand{\X}{{\Bbb X}}

\newcommand{\calC}{{\cal C}}
\newcommand{\calD}{{\cal D}}
\newcommand{\calH}{{\cal H}}
\newcommand{\calO}{{\cal O}}
\newcommand{\calR}{{\cal R}}
\newcommand{\calS}{{\cal S}}
\newcommand{\calU}{{\cal U}}

\newcommand{\frX}{{\frak X}}

\newcommand{\pair}[2]{\left\langle {#1}, \, {#2}\right\rangle}

\newcommand{\set}[2]{\left\{\left.#1\vphantom{#2}\:\right\vert\:#2\right\}}
\newcommand{\wt}{\widetilde}
\newcommand{\what}{\widehat}

\newcommand{\xx}{{\frak X}}

\newcommand{\lam}{{\lambda}}
\newcommand{\Lam}{{\Lambda}}
\newcommand{\ve}{{\varepsilon}}
\newcommand{\alp}{{\alpha}}  
\newcommand{\eps}{{\epsilon}} 

\newcommand{\CKXT}{\calC^\infty(K \backslash X_T)}   
\newcommand{\SKXT}{\calS(K \backslash X_T)}   

\newcommand{\slit}{\vspace{5mm}}
\newcommand{\mslit}{\vspace{3mm}}
\newcommand{\real}{{\rm Re}}

\newcommand{{\any}}{{}^\forall}
\newcommand{{\is}}{{}^\exists}
\newcommand{{\st}}{\; {\rm s.t.}\; }

\newcommand{\ol}[1]{\overline{#1}}

\newcommand{\hec}{{\calH(G,K)}}

\newcommand{\CKXU}{{\calC^\infty(K \backslash \frX_T/ U(T))}}   

\newcommand{\abs}[1]{\left\vert{#1}\right\vert}  

\newcommand{\dint}[1]{\displaystyle{\int_{#1}}}

\newcommand{\dprod}[1]{\displaystyle{\prod_{#1}}}

\newcommand{\dsqcup}[1]{\displaystyle{\bigsqcup_{#1}}}

\newcommand{\shita}[2]{\stackrel{\scriptstyle{#1}}{#2}}

\newcommand{\twomatrix}[4]{\begin{pmatrix}
                           {#1} & {#2}\\
                           {#3} & {#4}
                          \end{pmatrix}}

\newcommand{\twovector}[2]{\begin{pmatrix}{#1}\\{#2}\end{pmatrix}}

\newcommand{\twomatrixplus}[4]{\left(\begin{array}{c|c}
                           {#1}  & {#2}\\
                           \hline
                           {#3} &  {#4}
                          \end{array}\right)}
\newcommand{\twomatrixminus}[4]{\begin{array}{cc}
                           {#1}  & {#2}\\
                           {#3} &  {#4}
                          \end{array}}
\newcommand{\fourmatrix}[8]{\left(\begin{array}{cc|cc}
                           {#1} & {} & {#2} & {}\\
                           {} & {#3} & {} & {#4}\\
                           \hline
                           {#5} & {} & {#6} & {}\\
                           {} & {#7} & {} & {#8}
                          \end{array}\right)}

\newcommand{\proof}{{\hspace*{0.4cm} {\it Proof}.\ \enskip}}

\begin{document}

\title{Spherical functions on $U(2n)/\left(U(n) \times U(n)\right)$ and hermitian Siegel series}
\author{Yumiko Hironaka\\
{\small Department of Mathematics, }\\
{\small Faculty of Education and Integrated Sciences, Waseda University}\\
{\small Nishi-Waseda, Tokyo, 169-8050, JAPAN}
}

\date{}
\maketitle


\renewcommand{\thefootnote}{{}}
\makeatletter
\footnotetext{{}



E-mail: hironaka@waseda.jp

This research is partially supported by Grant-in-Aid for scientific Research (C):20540029. 
}
\makeatother

\setcounter{section}{-1}
\section{Introduction}
Let $\G$ be a reductive algebraic group and $\X$ an affine algebraic variety which is $\G$-homogeneous, where everything is assumed to be defined over a non-archimedian local field $k$ of characteristic $0$. We denote by $G$ and $X$ the sets of $k$-rational points of $\G$ and $\X$, respectively, take a maximal compact subgroup $K$ of $G$, and consider the Hecke algebra $\hec$. Then, a nonzero $K$-invariant function on $X$ is called {\it a spherical function on $X$} if it is an $\hec$-common eigenfunction.

Spherical functions on homogeneous spaces comprise an interesting topic to investigate and a basic tool to study harmonic analysis on $G$-space $X$. They have been studied as spherical vectors of distinguished models, Shalika functions and Whitakker-Shintani functions, there are close relation to the theory of automorphic forms, and spherical functions may appear in local factor of global object like Rankin-Selberg convolution and Eisenstein series. The theory of spherical functions has also an application of classical number theory, e.g. local densities of representations of quadratic forms or hermitian forms.  

To obtain explicit formulas of spherical functions is one of basic problems, and it has been done for  the group cases by I.~G.~Macdonald and afterwards by W.~Casselman by a representation theoretical method(\cite{Mac1}, \cite{Cas}). For homogeneous spaces, there are results mainly for   
the case that the space of spherical functions attached to each Satake parameter is of dimension one (e.g., \cite{CasS}, \cite{KMS}, \cite{Of}).  
The author gave general expressions of spherical functions on the basis of data of the group $G$ and functional equations of spherical functions when the dimension is not necessarily one, and a sufficient condition to have functional equations with respect to the Weyl group of $G$ (cf. \cite{French}). 

In the present paper, we investigate spherical functions on spaces $X_T$ for hermitian form $T$, where $X_T$ is a homogeneous space of the unitary group $\G = U(2n)$ and isomorphic to $U(2n)/U(n) \times U(n)$ over the algebraic closure of $k$ ($n$ is the size of $T$).
Here and henceforth we fix an unramified quadratic extension $k'$ of $k$ and consider hermitian forms and unitary groups with respect to the extension $k'/k$.  

In \S 1, we introduce the space $X_T = \frX_T/U(T)$ for each hermitian matrix $T$ of size $n$, and construct spherical functions $\omega_T(x; s) = \omega_T(x; z)$ on $X_T$, where $x \in X_T$ and $s, z \in \C^n$ are related by (\ref{change of var}).

We give the functional equations of $\omega_T(x; s)$ with respect to the Weyl group $W$
and determine the location of their possible poles and zeros (Theorem~\ref{th feq}, Theorem~\ref{th: W-inv}).
The Weyl group $W$ is isomorphic to $S_n \ltimes (\pm 1)^n$, and $S_n$ acts on the variable $z = (z_1, \ldots, z_n)$ by permutation of indices, and we may apply previous results on the spherical functions of hermitian forms to obtain the functional equations with respect to $S_n$. As for $\tau \in W$ corresponding to the remaining simple root, we need to consider the standard parabolic subgroup $P$ associated to $\tau$ and enlarge the space $\frX_T$ into $\wt{\frX_T}$ on which $P \times GL_1(k')$ acts. Different from the cases of the other simple roots, i.e., transpositions $(i\; i+1), \; 1 \leq i \leq n-1$, the functional equation with respect to $\tau$ does not come from that of a prehomogeneous vector space (cf. Remark~\ref{sect1 rem}, Theorem~\ref{th: tau}).  

Next we apply the general expression given in \cite{French} to the present case, and obtain the explicit formula for $\omega_T(x_T; z)$ for some diagonal $T$ and a particular point $x_T$ (Theorem~\ref{th: explicit}). Then, by sliding, we have the explicit formulas for general $T$ at many points (Theorem~\ref{th: many explicit}). 

In \S 4, we consider the spherical Fourier transform on the Schwartz space $\SKXT$, which is an integral transform employing the spherical function as kernel function, and show that the image is a free $\hec$-module of rank $2^{n-1}$. 

In \S 5, as an application, we consider hermitian Siegel series $b_\pi(T; t)$, relate them to our spherical functions $\omega_T(x; s)$. Then we give the `denominator part' of $b_\pi(T; t)$ and the functional equations of $b_\pi(T; t)$ by using the results in \S 2.  
A similar study for (symmetric) Siegel series has been done by F.~Sato and the author, but in that case we could not obtain the explicit formula by use of spherical functions. In the present case, we give the explicit functional by a specialization of functional equations of spherical functions $\omega_T(x; z)$. The existence of the functional equations was known in an abstract form as functional equations of Whitakker functions of $p$-adic groups by M.~L.~Karel, and explicit formulas have been given recently by T.~Ikeda (more precisely, see remarks in \S 5).

\vspace{2cm}
\section{Spherical function $\omega_T(x;s)$ on $\frX_T$ and $X_T$}

Let $k'$ be an unramified quadratic extension of a $p$-adic field $k$ with involution $*$, and 
for each $A = (a_{ij}) \in M_{mn}(k')$, we denote by $A^*$ the matrix $({a_{ji}}^*) \in M_{nm}(k')$. 
We fix a unit $\eps \in \calO_k^\times$ such that $k' = k(\sqrt{\eps})$ and $\ve - 1 \in 4\calO_k^\times$ (cf. \cite{Ome}, 63.3 and 63.4), and set 
\begin{eqnarray} \label{xi}
\xi = \frac{1 + \sqrt{\eps}}{2}.
\end{eqnarray}
Then $\{ 1, \; \xi \}$ forms an $\calO_k$-basis for $\calO_{k'}$, and $\set{\alp \in \calO_{k'}}{\alp^* = -\alp} = \sqrt{\eps} \calO_k$. 
We fix a prime element $\pi$ of $k$, and denote by $v_\pi(\; )$ the additive value on $k$, by $\abs{\; }$ the normalized absolute value on $k^\times$ with $\abs{\pi}^{-1} = q$ being the cardinality of the residue class field of $k$. 

We set
\begin{eqnarray*} \label{herm}
\calH_m = \set{A \in M_m(k')}{A^* = A}, \quad \calH_m^{nd} = \calH_m \cap GL_m(k').
\end{eqnarray*}
For $A \in \calH_m$ and $X \in M_{mn}(k')$, we write
$$
A[X] = X^* \cdot A = X^*AX \in \calH_n,
$$
and define the unitary group of $A$ by 
$$
U(A) = \set{g \in GL_m(k')}{A[g] = A}.
$$
In particular we set
$$
G = U(H_n) \quad \mbox{with } H_n = \twomatrix{0}{1_n}{1_n}{0}, 
\qquad U(m)  = U(1_m).
$$
For $T \in \calH_n^{nd}$, we set
\begin{eqnarray} \label{spaces}
&&
\xx_T = \set{x \in M_{2n,n}(k')}{H_n[x] = T}  \, \ni  x_T = \twovector{\xi T}{1_n},\nonumber\\
&&
 X_T = \xx_T \big{/}U(T). 
\end{eqnarray}
The group $G$ acts on $\xx_T$, as well as on $X_T$, through left multiplication, which is transitive by Witt's theorem for hermitian matrices (cf. \cite{Sch}, Ch.7, \S 9).

\begin{lem}
The stabilizer $G_0$ of $G$ at $x_T U(T) \in X_T$ is isomorphic to $U(T) \times U(T)$: 
\begin{eqnarray*} 
U(T) \times U(T) \stackrel{\sim}{\longrightarrow} G_0, \; (h_1, h_2) \longmapsto 
\wt{T}^{-1} 
\twomatrix {h_1^{* -1}}{0}{0}{h_2^{* -1}} \wt{T},
\end{eqnarray*}
where
$$
\wt{T} = 
\twomatrix{1_n}{\xi^* T}{1_n}{-\xi T}
\in GL_{2n}(k').
$$
\end{lem}

\proof
Since $\wt{T}x_T = \twovector{T}{0}$, we have, for any $h \in U(T)$,  
\begin{eqnarray} \label{h-part} 
x_Th = \wt{T}^{-1} \wt{T} x_T h = \wt{T}^{-1} \twovector{Th}{0} = \wt{T}^{-1} \twomatrix{h^{* -1}}{0}{0}{1} \wt{T}x_T.
\end{eqnarray}
Take any $g \in G$ such that $gx_T = x_T$. Then $\wt{T}g\wt{T}^{-1}$ stabilizes $\wt{T}x_T$ and belongs to $U(H_n[\wt{T}^{-1}])$, where 
$$
H_n[\wt{T}^{-1}] = H_n \left[ \twomatrix{\xi 1_n}{\xi^* 1_n}{ T^{-1} }{- T^{-1} } 
\right] =
\twomatrix{T^{-1}}{0}{0}{-T^{-1} },
$$
and we get
$$
\wt{T}g\wt{T}^{-1} =  \twomatrix{1_n}{0}{0}{d},  \mbox{ for some } d \in U(-T^{-1}). 
$$
Hence, together with (\ref{h-part}), we obtain the isomorphism stated as above.
\qed

\bigskip
We fix the Borel subgroup $B$ of $G$ as
\begin{eqnarray} \label{def Borel}
B = \set{\twomatrix{b}{0}{0}{b^{* -1}}\twomatrix{1_n}{a}{0}{1_n}}
  {\begin{array}{l}
      b \mbox{ is upper triangular of size }n, \\
      a + a^* = 0
     \end{array}},
\end{eqnarray}
ant the maximal compact sugbroup $K = G \cap GL_{2n}(\calO_{k'})$ of $K$, which satisfy $G = KB = BK$. 
We fix the $dk$ on $K$ and the left invariant Haar measure $dp$ on $B$ normalized by $\int_{K} dk = \int_{K \cap B} dp = 1$.
For each element $x \in \xx_T$, we denote by $x_2$ the lower half $n$ by $n$ block of $x$.
We define relative $B$-invariants on $\xx_T$ by
\begin{eqnarray} \label{rel inv}
f_{T, i}(x) = d_i(x_2 \cdot T^{-1}) = d_i (x_2 T^{-1} x_2^*), \quad 1 \leq i \leq n,
\end{eqnarray}
where $d_i(y)$ is the determinant of the upper left $i$ by $i$ block of a matrix $y$. It is easy to see, for $b \in B$,
\begin{eqnarray} \label{asso char}
f_{T,i}(bx) = \psi_i(b)f_{T,i}(x), \quad \psi_i(b) = N(d_i(b))^{-1},
\end{eqnarray}
where $N = N_{k'/k}$. 
Hence $f_{T,i}(x), \; 1 \leq i \leq n$ are relative $B$-invariants on $\frX_T$ associated with rational characters $\psi_i$ of $B$, and  we may regard them as relative $B$-invariants on $X_T$, since $f_{T, i}(xh) = f_{T, i}(x)$ for any $h \in U(T)$.
We set
\begin{eqnarray} \label{Xop}
\frX_T^{op} = \set{x \in X_T}{f_{T, i}(x) \ne 0, \; 1 \leq i \leq n}, \quad
X_T^{op} = \frX_T^{op}/U(T).
\end{eqnarray}

\begin{rem} \label{rem: realization}
{\rm
Though we may realize above objects as the sets of $k$-rational points of algebraic sets defined over $k$ and develop the arguments, we take down to earth way for simplicity of notations.
We only note here that $X_T^{op}$ (resp. $\frX_T^{op}$) becomes a Zariski open $B$-orbit in $X_T$ (resp. $B \times U(T)$-orbit in $\frX_T^{op}$) over the algebraic closure of $k$.
%
}
\end{rem}

\bigskip
We introduce a spherical function $\omega_T(x; s)$ on $\frX_T$ as well as on $X_T = \xx_T/U(T)$. 
For $x \in \frX_T$ and $s \in \C^n$, set 
\begin{eqnarray} \label{def sph f}
\omega_T(x; s) = \omega_T^{(n)}(x; s) = \dint{K}\, \abs{f_T(kx)}^{s+\ve} dk,
\end{eqnarray}
where $k$ runs over the set $\set{k \in K}{kx \in \frX_T^{op}}$,  
\begin{eqnarray}
&&
\ve = \ve_0 + (\frac{\pi\sqrt{-1}}{\log q},\ldots,\frac{\pi\sqrt{-1}}{\log q}), \quad \ve_0 = (-1, \ldots, -1, -\frac12) \in \C^n,\nonumber \\
&&
\label{def f(x)}
f_T(x) = \prod_{i=1}^n f_{T, i}(x), \qquad \abs{f_{T}(x)}^{s} = 
 \prod_{i=1}^n \abs{f_{T, i}(x)}^{s_i}. \nonumber
\end{eqnarray}
The right hand side of (\ref{def sph f}) is absolutely convergent if $\real(s_i) \geq -\real(\ve_i) = -\real(\ve_{0,i}), \; 1 \leq i\leq n$, and continued to a rational function of $q^{s_1}, \ldots, q^{s_n}$. 
We note here that 
$$
\abs{\psi(p)}^{\ve} \left( = \prod_{i=1}^n \abs{\psi_i(p)}^{\ve_i} \right) = \abs{\psi(p)}^{\ve_0} = \delta^{\frac12}(p),
$$ 
where $\delta$ is the modulus character on $B$ (i.e., $d(pp') = \delta(p')^{-1}dp$).

We denote by $\CKXT$ the space of left $K$-invariant functions on $X_T$, which can be identified with the space $\CKXU$ of left $K$-invariant right $U(T)$-invariant functions on $\frX_T$.  
The function $\omega_T(x; z)$ can be regarded as a function in $\CKXT$ and becomes a common eigenfunction by the action of the Hecke algebra $\calH(G, K)$ (cf. \cite{JMSJ} \S1, or \cite{French} \S 1).
In detail, 
the Hecke algebra $\hec$ is the commutative $\C$-algebra consisting of compactly supported two-sided $K$-invariant functions on $G$, acting on $\CKXT$ by the convolution product
\begin{eqnarray} \label{conv prod}
(\phi * \Psi) (x) = \int_{G}\, \phi(g) \Psi(g^{-1}x)dg, \quad (\phi \in \hec, \; \Psi \in \CKXT),
\end{eqnarray}
where $dg$ is the Haar measure on $G$ normalized by $\int_K dg = 1$, and we see 
\begin{eqnarray}  \label{eigen funt eq}
\left( \phi * \omega_T(\; ; s) \right)(x) = \lam_s(\phi) \omega_T(x; s),\quad (\phi \in \hec)
\end{eqnarray}
where $\lam_s$ is the $\C$-algebra homomorphism defined by
\begin{eqnarray*} \label{eigen funct}
&&
\lam_s : \hec \longrightarrow \C(q^{s_1}, \ldots, q^{s_n}), \\
&&
\quad
\phi \longmapsto \int_{B}\, \phi(p)\abs{\psi(p)}^{-s+\ve} dp, \quad \Big(\abs{\psi(p)}^{-s+\ve} = \prod_{i=1}^n \abs{\psi_i(p)}^{-s_i+\ve_i} \Big).
\end{eqnarray*}

\bigskip
We introduce a new variable $z$ which is related to $s$ by
\begin{eqnarray} \label{change of var}
s_i = -z_i + z_{i+1} \quad (1 \leq i \leq n-1), \quad
s_n = -z_n 
\end{eqnarray}
and write $\omega_T(x;z) = \omega_T(x; s)$. 
The Weyl group $W$ of $G$ relative to the maximal $k$-split torus in $B$ acts on rational characters of $B$ as usual (i.e., $\sigma(\psi)(b) = \psi(n_\sigma^{-1}b n_\sigma)$ by taking a representative $n_\sigma$ of $\sigma$), so $W$ acts on  $z \in \C^n$ and on $s \in \C^n$ as well.  We will determine the functional equations of $\omega_T(x; s)$ with respect to this Weyl group action.
The group $W$ is isomorphic to $S_n \ltimes C_2^n$, $S_n$ acts on $z$ by permutation of indices, and $W$ is generated by $S_n$ and $\tau: (z_1, \ldots, z_n) \longmapsto (z_1, \ldots, z_{n-1}, -z_n)$. Keeping the relation (\ref{change of var}), we also write $\lam_z(\phi) = \lam_s(\phi)$; then $\lam_z$ gives a $\C$-algebra isomorphism (the Satake isomorphism)
\begin{eqnarray} \label{satake iso}
\lam_z  &:& \hec \stackrel{\sim}{\longrightarrow} \C[q^{\pm 2z_1}, \ldots, q^{\pm 2z_n}]^W,\\
{} && \phi \longmapsto \dint{B}\, \phi(p) \prod_{i=1}^n\, \abs{N(p_i)}^{-z_i}\, \delta^{\frac12}(p) dp, \nonumber
\end{eqnarray}
where $p_i$ is the $i$-th diagonal component of $p \in B$, and the right hand side is the invariant subring of the Laurent polynomial ring by $W$.
\medskip
\begin{propos} \label{many sph}
Set $\calU = (\Z/2\Z)^{n-1}$ and 
$$
\wt{u} = (u_1\frac{\pi\sqrt{-1}}{\log q}, \ldots, u_{n-1}\frac{\pi\sqrt{-1}}{\log q}, 0) \in \C^n, \qquad u = (u_1, \ldots. u_{n-1}) \in \calU.
$$
Then 
$\omega_T(x; z + \wt{u}), \; u \in \calU$, are linearly independent for generic $z \in \C^n$ and correspond to the same eigenvalue through $\lam_z : \hec \longrightarrow \C$.
\end{propos}
 
\proof
The set $\frX_T^{op}$ is decomposed into the disjoint union of $B$-orbits as follows:
\begin{eqnarray*} \label{factor in frX-op}
&&
\frX_T^{op} = \dsqcup{u \in \calU}\, \frX_{T,u}, \\
&&
\frX_{T,u} = \set{x \in \frX_T^{op}}{v_\pi(f_{T,i}(x)) \equiv u_1 + \cdots + u_i\pmod{2}, \; 1 \leq i \leq n-1}.\nonumber
\end{eqnarray*}
We consider finer spherical functions 
$$
\omega_{T, u}({x}; s) = \dint{K}\, \abs{f_T(kx)}_u^{s+{\ve}}dk, \quad
\abs{f_T(y)}_u^{s+{\ve}} = \left\{ \begin{array}{ll} 
\abs{f_T(y)}^{s+{\ve}} & \mbox{if } {y} \in \frX_{T,u},\\
{} & {}\\
0 & \mbox{otherwise .}
\end{array}
\right.
$$
Then $\set{\omega_{T,u}(x, z)}{u \in \calU}$ are linearly independent for generic $z$ associated with the same $\lam_z$, where we keep the relation (\ref{change of var}) between $s$ and $z$. 
For each character $\chi$ of $\calU$, we may represent as follows
\begin{eqnarray} \label{omega(z_chi)}
\sum_{u \in \calU}\, \chi(u) \omega_{T,u}(x;s) = \omega_T(x; z_\chi),
\end{eqnarray}
where $z_\chi$ is obtained by adding $\frac{\pi\sqrt{-1}}{\log q}$ to $z_i$ for suitable $i$ according to $\chi$, and they are linearly independent (for generic $z$) as varying characters $\chi$. The result follows from this, since $\set{\omega_T(x; z_\chi)}{\chi \mbox{\, is a character of \,} \calU} = \set{\omega_T(x; z + \wt{u})}{u \in \calU}$. 
\qed

\bigskip
We note here the relation between $\omega_T(x; s)$ and $\omega_{T'}(y; s)$ when $T$ and $T'$ are equivalent under the action of $GL_n(k')$, which is easy to see. 

\begin{propos} \label{T[h]}
For $T \in \calH_n^{nd}$ and $h \in GL_n(k')$, we set $T' = T[h] \left( = h^*Th \right)$. Then
\begin{eqnarray*}  
&&
\frX_{T'} = \left( \frX_T \right)h, \quad X_{T'} = \frX_Th/U(T') \quad \mbox{and} \quad
f_{T', i}(xh) = f_{T, i}(x) \quad (x \in \xx_T), \nonumber
\end{eqnarray*}
and 
\begin{eqnarray*}  
&&
{\omega}_{T'}(xh; s) = {\omega}_T(x; s), \qquad (x \in \xx_T).
\end{eqnarray*}
\end{propos}


\bigskip
By using a result on spherical functions on the space of hermitian forms, we obtain the following theorem.

\begin{thm} \label{aim S_n} 
Set 
\begin{eqnarray} \label{G1(z)}
G_1(z) = \dprod{1 \leq i < j \leq n}\, \frac{q^{z_j} + q^{z_i}}{q^{z_j} - q^{z_i-1}} .
\end{eqnarray}
Then, for any $T \in \calH_n^{nd}$, the function $G_1(z) \cdot \omega_T(x; z)$ is invariant under the action of $S_n$ on $z$.  
\end{thm}

\proof
By the embedding 
\begin{eqnarray*} \label{embed}
K_0 = GL_n(\calO_{k'})
\longrightarrow K, \quad h \longmapsto \wt{h} = \twomatrix{h^{* -1}}{0}{0}{h}, 
\end{eqnarray*}
and the normalized Haar measure $dh$ on $K_0$, we obtain, for $s \in \C^n$ satisfying $\real(s_i) \geq -\real(\ve_i), \; 1 \leq i \leq n$,
\begin{eqnarray*}
\omega_T(x;z) &=& \omega_T(x;s) = 
\dint{K_0}\,dh \dint{K}\, \abs{f_T(kx)}^{s+\ve} dk
\nonumber \\
&=&
\dint{K_0}\,dh \dint{K}\, \abs{f_T(\wt{h}kx)}^{s +\ve} dk
=
\dint{K}\, \dint{K_0}\, \abs{f_T(\wt{h}kx)}^{s + \ve} dh dk\nonumber\\
&=&
\dint{K}\, \zeta^{(n)}(D(kx); s) dk.
\label{id with herm}
\end{eqnarray*}
Here $D(kx) = (kx)_2 \cdot T^{-1} \in \calH_n^{nd}$ for $\set{k \in K}{kx \in \frX_T^{op}}$, and $\zeta^{(n)}(y; s)$ is a spherical function on $\calH_n^{nd}$ defined by
\begin{eqnarray*}
\zeta^{(n)}(y; s) = \dint{K_0}\, \prod_{i=1}^n\, \abs{d_i(h\cdot y)}^{s_i+\ve_i} dh, \quad (h \cdot y = hyh^*),
\end{eqnarray*}
where $h$ runs over the set $\set{h \in K_0}{d_i(h\cdot y) \ne 0, \; 1 \leq i \leq n}$. Keeping the relation between $s$ and $z$ as before, the assertion of Theorem~\ref{aim S_n} follows from the next proposition. 
\qed
%
\begin{propos} {\rm (cf. {\rm \cite{Y-III} or \cite{Hamb}})} \label{herm}
For any $y \in \calH_n^{nd}$, the function $G_1(z) \cdot \zeta^{(n)}(y; s)$ 
is holomorphic for $z \in \C^n$ and invariant under the action of $S_n$, where $G_1(z)$ is defined as in (\ref{G1(z)}).
\end{propos}

In \cite{Hamb} \S 4.2, we considered a modified function
\begin{eqnarray*} \label{Hamb}
\omega^{(H)}(y; s) = \dint{K_0}\, \prod_{i=1}^n\, \chi_\pi(d_i(h \cdot y)) \abs{d_i(h\cdot y)}^{s_i+\ve'_i} dh,
\end{eqnarray*}
where $\chi_\pi(a) = (-1)^{v_\pi(a)}$ for $a \in k^\times$ and $\ve' = (-1, \ldots, -1, \frac{n-1}2)$. The function $\zeta^{(n)}(x; s)$ satisfies the same functional properties as $\omega^{(H)}(y; s)$, since $\omega^{(H)}(y; s) = \abs{\det(y)}^{\frac{n}{2}} \zeta^{(n)}(y; s)$. 

\begin{rem} \label{sect1 rem}
{\rm 
For the transposition  $\tau_i = (i\; i+1) \in W$, \; $1 \leq i \leq n-1$, the following functional equations hold by Theorem~\ref{aim S_n}
\begin{eqnarray} \label{feq of tau_i}
\omega_T(x; z) = 
\frac{ 1 - q^{z_i -z_{i+1}-1} }{ q^{z_i-z_{i+1}} - q^{-1} } \times \omega_T(x; \tau_i(z)), \quad 1 \leq i \leq n-1.
\end{eqnarray} 
On the other hand, one may obtain (\ref{feq of tau_i}) directly in the similar way to the case of $\tau$ in \S~2, where the sufficient condition in \cite{French}-\S 3 for having a functional equation with respect to $\tau_i$ is satisfied and the  Gamma factor in (\ref{feq of tau_i}) is essentially the same to that of the zeta function of prehomogeneous vector space $(U \times GL_1(k'), (k')^2)$, where $U \cong U(2)$ or $U(Diag(1, \pi))$.  
Then Theorem~\ref{aim S_n} follows from (\ref{feq of tau_i}), through the similar line to the proof of Proposition~\ref{herm}. In fact, Proposition~\ref{herm} was proved by using functional equations of type (\ref{feq of tau_i}).
}
\end{rem}

\vspace{2cm}

\setcounter{section}{1}
\setcounter{equation}{0}
\section{Functional equations, possible zeros and poles}
We calculate the functional equation for $\tau \in W$, and give the functional equations with respect to the whole $W$.

\bigskip
\noindent
{\bf 2.1.} First we calculate the spherical function for $n = 1$.
We note the data for $n =1$, which will be used also in \S 2.2.
\begin{eqnarray}
&&
G = U(1,1) = 
\set{\twomatrix{\alp}{0}{0}{\alp^{* -1}} \twomatrix{1}{w\sqrt{\eps}}{v\sqrt{\eps}}{1+vw\eps}}{\alp \in {k'}^\times, \; v, w \in k}  \nonumber \\
&&
\hspace*{2.7cm} \cup
\set{\twomatrix{\alp}{0}{0}{\alp^{* -1}} \twomatrix{0}{1}{1}{v\sqrt{\eps}}}{\alp \in {k'}^\times, \; v \in k},\nonumber \\
&&
K = K_1 = K_{1,1} \cup K_{1,2}, \mbox{ where} \nonumber \\
&&
\hspace*{.7cm}
K_{1,1} = \set{\twomatrix{\alp}{0}{0}{\alp^{*-1}}\twomatrix{1}{v/\sqrt{\eps} } {u\sqrt{\eps}}{1+uv} }{ \alp \in \calO_{k'}^\times, \; u, v \in \calO_k}, \nonumber \\
&&
\hspace*{.7cm}
K_{1,2} = \set{ \twomatrix{\alp}{0}{0}{\alp^{*-1}}\twomatrix{\pi u \sqrt{\eps}}{1+\pi uv}{1}{v/\sqrt{\eps} } }
{ \alp \in \calO_{k'}^\times, \; u, v \in \calO_k}, \label{K1}
\end{eqnarray}
and
\begin{eqnarray*}
\omega^{(1)}_T(x; s) = \dint{K_1}\, \chi_\pi(f_1(hx)) \abs{f_1(hx)}^{s-\frac12} dh,
\end{eqnarray*}
where 
$f_1(x) = \det(T)^{-1} N(x_2)$ for $x \in \xx_T$, 
$\chi_\pi(a) = (-1)^{v_\pi(a)}$ for $a \in k^\times$ and $dh$ is the Haar measure on $K_1$.

\begin{propos} \label{size1 prop}

\noindent
{\rm (i)} The set
\begin{eqnarray*} \label{rep n=1}
\set{x_e = \twovector{\pi^e}{\xi \pi^{t - e}}}{e \in \Z, \; 2e \leq t }, \quad \left(\xi = \frac{1 + \sqrt{\eps}}{2} \right)
\end{eqnarray*}
forms a complete set of representatives of $K_1 \backslash \xx_T$ for $T = \pi^t$. 

\noindent
{\rm (ii)}
For $x_e \in \xx_T$ with $T = \pi^t$ as above, one has
\begin{eqnarray*}
\omega^{(1)}_T(x_e; s) &=& 
\frac{(-1)^t q^{e-\frac12{t} } }{1+q^{-1}}
\times
 \frac{q^{(t-2e+1)s}(1-q^{-2s-1}) - q^{-(t-2e+1)s}(1-q^{2s-1})}{q^s - q^{-s}} .
\end{eqnarray*}

\noindent
{\rm (iii)}
For any $T \in \calH_1^{nd}$, $\omega^{(1)}_T(x; s)$ is holomorphic for all $s \in \C$ and satisfies the functional equation
$$
\omega^{(1)}_T(x; s) = \omega^{(1)}_T(x; -s).
$$
\end{propos}

\proof
We recall that $\{1, \xi \}$ forms an $\calO_k$-basis of $\calO_{k'}$ and $Tr_{k'/k}(\xi) = 1$.  
Multiplying a suitable element in $K_1$ of type
$$
\twomatrix{1}{0}{u\sqrt{\ve}}{1}\twomatrix{\alp}{0}{0}{\alp^{* -1}} \quad \mbox{or} \quad
\twomatrix{1}{0}{u\sqrt{\ve}}{1}\twomatrix{0}{\alp}{\alp^{* -1}}{0} \quad (u \in \calO_k, \; \alp \in \calO_{k'}^\times)
$$
one may make any $x \in \frX_T$ into some $x_e$ in the given set, 
and the explicit formula in (ii) shows there is no redundancy within it.  

Take $x_e$ as above. For $h \in K_{1,1}$ written as in (\ref{K1}), since we have 
$$
(hx_e)_2 = \alp^{* -1}(u\sqrt{\eps}\pi^e + (1+uv)\xi\pi^{t-2e})) 
= \alp^{* -1} \pi^e\left( -u + ((1+uv)\pi^{t-2e} + 2u)\xi \right),
$$
we see 
$$
v_\pi(f_1(kx_e)) = -t+2e + 2\min\{v_\pi(u), \, t-2e \}.
$$
Since $vol(K_{1,1}) = (1 + q^{-1})^{-1}$, we obtain
\begin{eqnarray*}
\lefteqn{ \dint{K_{1,1}}\, \chi_\pi(f_1(hx_e)) \abs{f_1(hx_e)}^{s-\frac12} dh}\\
& = &
\frac{(-1)^t q^{(t-2e)(s-\frac12)}}{1+q^{-1}}\cdot  
\sum_{r \geq 0}q^{-r}(1-q^{-1}) q^{-2\min\{r,\, t-2e \}(s-\frac12)} \\
&=&
\frac{(-1)^t q^{(t-2e)(s-\frac12)}}{1+q^{-1}} \cdot \left( 
\frac{(1-q^{-1})(1-q^{-2(t-2e)s})}{1-q^{-2s}} + q^{-2(t-2e)s} \right).
\end{eqnarray*}
For $h \in K_{1,2}$ written as in (\ref{K1}), since we have
$$
(hx_e)_2 = \alp^{* -1} ( \pi^e + v/\sqrt{\ve} \pi^{t-e}\xi)
= (\alp^*\sqrt{\ve})^{-1} \pi^e \left( -1 + (2+v\pi^{t-2e}) \xi \right),
$$
we see
$v_\pi(f_1(hx_e)) = -t+2e$ and 
$$
\dint{K_{1,2}}\, \chi_\pi(f_1(hx_e)) \abs{f_1(hx_e)}^{s-\frac12} dh = \frac{q^{-1}}{1+q^{-1}} \cdot (-1)^t q^{(t-2e)(s-\frac12)}.
$$
Thus we obtain
\begin{eqnarray*}
\lefteqn{\omega^{(1)}_T(x_e; s) = \frac{(-1)^t q^{(t-2e)(s-\frac12)}}{1+q^{-1}} \frac{1}{1-q^{-2s}} \cdot \left( 1 - q^{-2s-1} + q^{-2(t-2e)s-1} - q^{-2(t-2e+1)s)}
 \right) }\\
 &=&
\frac{(-1)^t q^{e-\frac12{t} } }{1+q^{-1}} \cdot \frac{1}{q^s - q^{-s}} \cdot
\left( q^{(t-2e+1)s}(1-q^{-2s-1}) - q^{-(t-2e+1)s}(1-q^{2s-1}) \right), 
\end{eqnarray*}
which proves (ii) and (iii) for $T = \pi^t$. Then, by Proposition~\ref{T[h]} we obtain the assertion (iii) for general $T \in \calH_1^{nd}$, since $N(\calO_{k'}^\times) = \calO_k^\times$.
\qed

\begin{rem}\label{rem size1}
{\rm 
In $z$-variable, the assertion in Proposition~\ref{size1 prop} becomes as follows, where $z = -s$ and $W = \{1, \tau \}$. For $T = \pi^t$,
\begin{eqnarray*}
&&
\omega^{(1)}_T(x_e; z) =  
\frac{(-1)^t q^{e-\frac12{t} } }{1+q^{-1}} \times \left\{ \frac{q^{-(t-2e)z}(1-q^{2z-1})}{1-q^{2z}} +\frac{q^{(t-2e)z}(1-q^{-2z-1})}{1-q^{-2z}}  \right\}, \quad(2e \leq t);
\end{eqnarray*}
and for any $T \in \calH_1^{nd}$ and $x \in \frX_T$,
\begin{eqnarray*}
\omega^{(1)}_T(x; z) = \omega^{(1)}_T(x; \tau(z)).  
\end{eqnarray*}
}
\end{rem}

\bigskip
\noindent
{\bf 2.2.} 
Assume that $n \geq 2$ and 
set
\begin{eqnarray*}
&&
w_\tau = \fourmatrix{1_{n-1}}{}{0}{1}{}{1_{n-1}}{1}{0} \in G, 
\end{eqnarray*}
then $w_\tau$ gives the element $\tau \in W$ such that $\tau(z) = (z_1, \ldots, z_{n-1}, -z_n)$.
The main purpose of this subsection is to prove the following.

\begin{thm} \label{th: tau}
For any $T \in \calH_n^{nd}$, the spherical function satisfies 
$$
\omega_T(x; z) = \omega_T(x; \tau(z)).
$$
\end{thm}

\bigskip
The standard parabolic subgroup $P$ attached to $\tau$, in the sense of \cite{Borel} \S 21.11, is given as follows:
\begin{eqnarray} 
\lefteqn{P =  B \cup B w_\tau B } \nonumber \\
&=& \label{P-tau}
\left\{ \left. 
\fourmatrix {q} {} {a} {b} {} {q^{*-1}} {c} {d} 
 \twomatrixplus
 {\twomatrixminus{1_{n-1}}{\alp}{}{1}} {} {}
 {\twomatrixminus{1_{n-1}}{}{-\alp^*}{1}}
 \twomatrixplus
 {1_n} { \twomatrixminus{\gamma}{\beta}{-\beta^*}{0} } {} {1_n} \in G 
\right\vert \right. \nonumber \\
&& \hspace*{4cm} 
\left. \vphantom{
 \twomatrixplus
 {\twomatrixminus{1_{n-1}}{\alp}{}{1}} {} {}
 {\twomatrixminus{1_{n-1}}{}{-\alp^*}{1}}
} 
\begin{array}{l} 
q \mbox{ is upper triangular in }GL_{n-1}(k'),\\
\twomatrix{a}{b}{c}{d} \in U(H_1), \;
 \alpha, \beta \in M_{n-1,1}(k'),\\
 \gamma \in M_{n-1}(k'), \; \gamma + \gamma^* = 0
  \end{array}
 \right\},  
\end{eqnarray}
where each empty place in the above expression means zero-entry.

Since it suffices to show Theorem~\ref{th: tau} for diagonal $T$'s (cf. Proposition~\ref{T[h]}),
we fix a diagonal $T \in \calH_n^{nd}$ and write $f_i(x) = f_{T,i}(x)$ for simplicity of notations.
We consider the following action of $\wt{P} = P \times GL_1$ on $\wt{\frX}_T = \frX_T \times V$ with $V = M_{21}(k')$:
$$
(p,r)\star (x,v) = (px, \rho(p)vr^{-1}), \qquad (p,r) \in \wt{P}, \; (x,v) \in \wt{\frX}_T,
$$
where $\rho(p) = \twomatrix{a}{b}{c}{d}$ for the decomposition of $p \in P$ as in (\ref{P-tau}). 
We define
\begin{eqnarray}  \label{g(x,v)}
g(x, v) &=& \det\left[ \twomatrixplus{1_{n-1}}{0}{0}{{}^tv}\twovector{x_2}{-y} \cdot T^{-1}\right],
\quad (x, v) \in \wt{\frX}_T,
\end{eqnarray}
where the first matrix in the right hand side is of size $(n, n+1)$, $x_2$ is the lower half $n$ by $n$ block of $x$ (the same as before) and  
$y$ is the $n$-th row of $x$.

\begin{lem}   \label{lem g(x,v)}
Let $g(x,v)$ be the function on $\wt{\frX}_T = \frX_T \times V$ defined by (\ref{g(x,v)}). 

\noindent
{\rm (i)} $g(x,v)$ is a relative $\wt{P}$-invariant on $\wt{\frX}_T$ associated with character $\wt{\psi}$: 
$$
\wt{\psi}(p,r) = \psi_{n-1}(p) N(r)^{-1}, \quad (p,r) \in \wt{P} = P \times GL_1,
$$
where $\psi_{n-1}$ is well-defined on $P$, and satisfies 
$$
g(x, v_0) = f_n(x), \qquad v_0 = \twovector{1}{0} \in V
$$
\noindent
{\rm (ii)} $g(x,v)$ is expressed as  
\begin{eqnarray} \label{g-D(x)}
g(x,v) = D(x)[v],
\end{eqnarray}
with some hermitian matrix 
\begin{eqnarray} \label{D(x)}
D(x) = \twomatrix{a(x)}{\beta(x)}{\beta(x)^*}{d(x)} \qquad 
(a(x), d(x) \in k, \; \beta(x) \in k'),
\end{eqnarray}
such that $\det(D(x)) = 0$ and $\trace(\beta(x)) = -f_{n-1}(x)$ for $x \in \frX_T$, where $\trace$ is the trace $\trace_{k'/k}$.
\end{lem}

\proof
(i) It is easy to see that $g((1,r) \star (x,v)) = N(r)^{-1} g(x,v)$. In order to examine the action of $P$, we write an element $p \in P$ and $x \in \frX_T$ as follows 
\begin{eqnarray*}
&&
p = 
\left( \begin{array}{cc|cc}
q \, & \cdot & \cdot\,  & \cdot\\[1mm]
0 \, & a & t \, & b\\
\hline
0 & 0 & q^{* -1} & 0\\[1mm]
0 & c & \mu & d
\end{array} \right), \qquad x = \begin{pmatrix}\cdot\\y\\[0.1mm] x'\\z 
\end{pmatrix}, 
\end{eqnarray*}
where $q \in GL_{n-1}, \; t, \mu \in M_{1, n-1}, \; 
\twomatrix{a}{b}{c}{d} \in U(H_1)$,  $x' \in M_{n-1,n}$, and $y, z \in M_{1,n}$. 
Then we obtain
\begin{eqnarray*}
g((p,1)\star(x,v)) &=& 
\det\left[ \twomatrixplus{1_{n-1}}{0}{0}{{}^tv \twomatrix{a}{c}{b}{d}}
 \left( \begin{array}{c} {q^{* -1}x'}\\ {cy+\mu x'+dz}\\ -(ay+t x' + bz) 
  \end{array}\right)
 \cdot T^{-1} \right] \\[2mm]
&=&
\det\left[ \twomatrixplus{1_{n-1}}{0}{0}{{}^tv\twomatrix{a}{c}{b}{d}}
\twomatrixplus{q^{* -1}}{0}{\begin{array}{c}\mu\\ -t \end{array}}{\twomatrixminus{d}{-c}{-b}{a}}
\left( \begin{array}{c}
x' \\ z \\ -y \end{array} \right) \cdot T^{-1} \right]\\[2mm]
&=&
\det\left[ \twomatrixplus {q^{* -1}}{0} 
              { {}^tv\twomatrix{a}{c}{b}{d}\twovector{\mu}{-t} } {\ve{}^tv}
\left( \begin{array}{c}
x' \\ z \\ -y \end{array} \right) \cdot T^{-1} \right]\quad(\ve = ad-bc \in \calO_{k'}^1) \\[2mm]
&=&
\det\left[ \twomatrixplus{q^{* -1}}{0}{{}^tv\twomatrix{a}{c}{b}{d}\twovector{\mu}{-t}}{\ve} 
\twomatrixplus{1_{n-1}}{0}{0}{{}^tv}
\left( \begin{array}{c}
x' \\ z \\ -y \end{array} \right) \cdot T^{-1} \right] \\[2mm]
&=&
N(\det(q))^{-1} g(x,v) = \psi_{n-1}(p) g(x,v).
\end{eqnarray*}
Hence $g(x,v)$ is a relative $\wt{P}$-invariant on $\wt{\frX}_T$ associated with character $\wt{\psi}$.

(ii)
Since $g(x,v)$ is a linear form with respect to both $v_1, v_2$ and $v_1^*, v_2^*$, and $g(x,v)^* = g(x,v)$, we have an expression (\ref{g-D(x)}) with some $D(x) \in \calH_2$. 
Writing $T = Diag(t_1, \ldots,t_n)$, we have
\begin{eqnarray}
g(x_T, v) &=& (t_1\cdots t_n)^{-1} (v_1- \xi t_n v_2)(v_1^*- \xi^* t_nv_2^* ) \nonumber \\[1mm]
&=& \label{g(xT)}
\twomatrix{(t_1\cdots t_n)^{-1}} {-\xi(t_1\cdots t_{n-1})^{-1}}
{-\xi^*(t_1\cdots t_{n-1})^{-1}} {N(\xi) (t_1\cdots t_{n-1})^{-1}t_n}[v],
 \nonumber \\
&=&
\twomatrix{(t_1\cdots t_n)^{-1}} {-\xi f_{n-1}(x_T)}
{-\xi^* f_{n-1}(x_T)} {N(\xi) (t_1\cdots t_{n-1})^{-1}t_n}[v],
\label{D(xT)}
\end{eqnarray}
in particular $\det(D(x_T)) = 0$. %
Since $g(x,v)$ is a relative $\wt{P}$-invariant associated with $\wt{\psi}$ by (i), we see
\begin{eqnarray} \label{D(px)}
D(px) = \psi(p)D(x)[\rho(p)^{-1}], \quad (p \in P)
\end{eqnarray}
and we have
\begin{eqnarray*} \label{for P}
\det(D(px_T)) = 0.
\end{eqnarray*}
Since $\frX_T^{op}$ is a $B$-orbit over the algebraic closure of $k$(cf. Remark~1.2), we have
\begin{eqnarray*}
\det(D(x)) = 0, \quad \mbox{ for any }x \in \frX_T^{op}.
\end{eqnarray*}
For an element $b$ of $B$, $\rho(b)$ can be written as follows (cf. (\ref{def Borel}))
$$
\rho(b) = \twomatrix{\gamma}{\gamma u\sqrt{\eps}}{0}{\gamma^{* -1}}, \quad \gamma \in {k'}^\times, \; u \in k,
$$
and when we express $D(x)$ and $D(bx)$ as in (\ref{D(x)}), we have by (\ref{D(px)})
\begin{eqnarray*}
\beta(bx) = \psi_{n-1}(b)(-a(x)u\sqrt{\ve} + \beta(x)),
\end{eqnarray*}
hence $\trace(\beta(bx)) = \psi_{n-1}(b) \trace({\beta(x)})$ and $\trace(\beta(x_T)) = -f_{n-1}(x_T)$
 by (\ref{g(xT)}). 
Thus $\trace(\beta(x)) = - f_{n-1}(x)$ for $x \in \frX_T^{op}$.
\qed

\bigskip
For $A \in \calH_2$ and $s \in \C$, we define
\begin{eqnarray*} \label{def zeta K1}
\zeta_{K_1}(A; s) = \dint{K_1}\, \abs{d_1(h \cdot A)}^{s-\frac12} dh,
\end{eqnarray*}
where $h\cdot A = hAh^*$ and $dh$ is the normalized Haar measure on $K_1 = U(H_1)\cap GL_2(\calO_{k'})$, which is absolutely convergent if $\real(s) \geq \frac12$.

\begin{lem} \label{prop : zeta K1}
Assume $x \in \frX_T^{op}$ and $D(x)$ is given by (\ref{g-D(x)}).  Set $m = \min\{v_\pi(a(x)), \; v_\pi(d(x)) \}$ and $t = v_\pi(\beta(x)) - m$ for the expression of $D(x)$ as in (\ref{D(x)}).
Then $t \geq 0$ and 
$$
\zeta_{K_1}(D(x); s) = \frac{q^{\frac{m}{2}}}{1+q^{-1}} \cdot \abs{f_{n-1}(x)}^s
\cdot \frac{q^{(t+1)s}(1 - q^{-2s-1}) - q^{-(t+1)s}(1 - q^{2s-1})}{q^s - q^{-s}}.
$$
In particular, one has the functional equation 
\begin{eqnarray} \label{feq zetaK1}
\zeta_{K_1}(D(x); s) = \abs{f_{n-1}(x)}^{2s} \cdot \zeta_{K_1}(D(x); -s).
\end{eqnarray}
\end{lem}

\proof
Take an $x \in \frX_T^{op}$, write $D(x)$ as in (\ref{D(x)}), and set $m$ as above.
Then $\beta(x)$ can be written as 
$$
\beta(x) = b_1 + \xi b_2, \quad b_1, b_2 \in k, \; m \leq \min\{b_1, b_2 \}, \; \trace(\beta(x)) = 2b_1 + b_2 = -f_{n-1}(x).
$$
Then, by the action of $K_1$ on  $\calH_2$, we see $D(x)$ is $K_1$-equivalent to 
\begin{eqnarray} \label{gen-k} 
\pi^m \twomatrix{1}{\xi b}{\xi^* b}{N(\xi) b^2}, \quad  b = \pi^{-m}\trace(\beta(x)) \in \calO_k,
\end{eqnarray}
and if $k$ is nondyadic, it is $K_1$-equivalent to 
\begin{eqnarray} \label{odd-k} 
\pi^m \twomatrix{1}{\frac{1}{2}b}{\frac{1}{2}b}{\frac14 b^2},\quad b = \pi^{-m}\trace(\beta(x)) \in \calO_k.
\end{eqnarray}
We denote by $A$ the matrix given in (\ref{gen-k}) (resp. in (\ref{odd-k}) ) if $k$ is dyadic (resp. nondyadic), then $\zeta_{K_1}(D(x); s) = \zeta_{K_1}(A; s)$.  
We recall the data for $K_1 = K_{1,1} \cup K_{1,2}$ in (\ref{K1}).

For $h = \twomatrix{\alp}{1}{\alp^{* -1}}{1} \twomatrix{1}{v/\sqrt{\eps} } 
{u\sqrt{\eps}}{(1+uv)} \in K_{1,1}$, we have
\begin{eqnarray*}
d_1(h \cdot A) =  \left\{\begin{array}{ll}
\pi^m N(\alp) (1 - \dfrac{b^2v^2}{4\eps}) & \mbox{if $k$ is nondyadic},\\[2mm]
\dfrac{\pi^m N(\alp)}{\eps} (\eps - \eps bv - N(\xi)b^2v^2 ) & \mbox{if $k$ is dyadic},
\end{array} \right.
\end{eqnarray*}
and $v_\pi(d_1(h\cdot A)) = m$ for any $h \in K_{1,1}$,
where we recall that $\eps \in 1 + 4\calO_k^\times$.

For $h = \twomatrix{\alp}{1}{\alp^{* -1}}{1} \twomatrix{\pi u \sqrt{\eps}}{1 + \pi u v}{1}{v/\sqrt{\eps} } \in K_{1,2}$ (cf. (2.1)), we have
\begin{eqnarray*}
d_1(h \cdot A) = \left\{  \begin{array}{ll}
\pi^m N(\alp) (-\eps \pi^2u^2  + (1+\pi uv)^2 b^2/4) & \mbox{if $k$ is non dyadic},\\[2mm]
\pi^m N(\alp) (-\eps \pi^2u^2 + (1+\pi uv)\eps \pi ub + (1+\pi uv)^2N(\xi) b^2) & \mbox{if $k$ is dyadic},
\end{array} \right.
\end{eqnarray*}
and $v_\pi(d_1(h\cdot A)) = m + 2\min\{v_\pi(b), v_\pi(u)+1 \}$. 

Set $t = v_\pi(b)$. If $t = 0$, it is clear that $\zeta_{K_1}(A; s) = q^{-m(s-\frac12)}$.  
If $t > 0$, then we obtain
\begin{eqnarray*}
\zeta_{K_1}(A; s) &=& 
\frac{1}{1+q^{-1}} a^{-m(s-\frac12)} + 
\frac{q^{-1}}{1 + q^{-1}} \left( \sum_{\ell=0}^{t-1}\, q^{-\ell}(1 - q^{-1})q^{-(m+2+2\ell)(s-\frac12)}  + q^{-t}q^{-(m+2t)(s-\frac12)}\right)
 \nonumber \\
&=&
\frac{q^{-m(s-\frac12)}}{(1 + q^{-1})} \times
\left( 1 + \frac{q^{-2s} - q^{-2s-1} + q^{-2t s - 1} - q^{-2(t +1)s}}   {1 - q^{-2s}} \right) \\
&=&
\frac{q^{-(m+t)s + \frac{m}{2}}}{(1 + q^{-1})} \times 
\frac{q^{(t+1)s}(1-q^{-2s-1}) - q^{-(t+1)s}(1 - q^{2s-1})}{q^s - q^{-s}},
\end{eqnarray*}
and the latter two expressions are valid also for $t = 0$.
Since $\pi^m b = \trace{\beta(x)} = -f_{n-1}(x)$, we have
\begin{eqnarray*}
\zeta_{K_1}(D(x); s) &=& \zeta_{K_1}(A; s)\\
&=&
\frac{ q^{\frac{m}{2}}}{1 + q^{-1}} \abs{f_{n-1}(x)}^s \times
\frac{q^{(t+1)s}(1-q^{-2s-1}) - q^{-(t+1)s}(1-q^{2s-1})}{q^s - q^{-s}}.
\end{eqnarray*}
The identity (\ref{feq zetaK1}) follows from the above explicit formula.
\qed

\bigskip
Now we will prove Theorem~\ref{th: tau}. We consider the embedding
\begin{eqnarray*}
K_1 \longrightarrow K = K_n, \quad h = \twomatrix{a}{b}{c}{d} \longmapsto 
\wt{h}= \fourmatrix {1_{n-1}} {} {a} {b} {} {1_{n-1}} {c} {d} .
\end{eqnarray*}
Then we have
\begin{eqnarray*}
\omega_T(x; s) &=& \dint{K_1}\, dh \dint{K} \abs{f(kx)}^{s+\ve} dk\\
&=&
\dint{K_1}dh \dint{K} \abs{f(\wt{h}kx)}^{s+\ve}dk\\
&=&
\dint{K}\, \chi_\pi(\dprod{i < n}f_i(kx)) \dprod{i<n}\,\abs{f_i(kx)}^{s_i-1} 
\left( \dint{K_1} \chi_\pi(f_n(\wt{h}kx))\abs{f_n(\wt{h}kx)}^{s_n-\frac12} dh \right) dk.
\end{eqnarray*}

By definition of $f_n(x)$ and $g(x,v)$ and Lemma~\ref{lem g(x,v)}, we have 
for $h = \twomatrix{a}{b}{c}{d} \in K_1
$\begin{eqnarray*}
f_n(\wt{h}x) &=& 
\det\left[  \twovector{x'}{cy+dz} \cdot T^{-1}  \right] 
= g(x, \twovector{d}{-c}) \\
&=& (d^*\; -c^*) D(x)\twovector{d}{-c}
= d_1(h^{* -1} \cdot D(x)).
\end{eqnarray*}
Since $\set{h^{* -1}}{h \in K_1} = K_1$, we have
\begin{eqnarray*}
\omega_T(x; s) = \dint{K}\,\chi_\pi(\dprod{i < n}f_i(kx)) \dprod{i<n}\,\abs{f_i(kx)}^{s_i-1} \zeta_{K_1}(D(kx); s_n + \textstyle{\frac{\pi\sqrt{-1}}{\log q}}) dk,  
\end{eqnarray*}
and by Lemma~\ref{prop : zeta K1}, we obtain
\begin{eqnarray*}
\lefteqn{\omega_T(x; s)}\\ 
&=&
\dint{K}\, \chi_\pi(\dprod{i < n}f_i(kx)) \dprod{i\leq n-2}\,\abs{f_i(kx)}^{s_i-1} 
\cdot 
\abs{f_{n-1}(kx)}^{s_{n-1} + 2s_n -1} \\
&& \times \zeta_{K_1}(D(kx); -s_n + \frac{\pi\sqrt{-1}}{\log q}) dk\\
&=&
\omega_T(x; s_1, \ldots, s_{n-2}, s_{n-1}+2s_n, -s_n). 
\end{eqnarray*}
In variable $z$, we have
\begin{eqnarray*}
\omega_T(x; z) = \omega_T(x; \tau(z)), \qquad \tau(z) = (z_1,\ldots,z_{n-1}, -z_{n}),
\end{eqnarray*}
which completes the proof.
\qed

\bigskip
\noindent
{\bf 2.3.}
In order to describe functional equations of $\omega_T(x; z)$ with respect to $W$, we prepare some notations.
We denote by $\Sigma$ the set of roots of $G$ with respect to the $k$-split torus of $G$ contained in $B$ and by $\Sigma^+$ the set of positive roots with respect to $B$. We may understand $\Sigma$ as a subset in $\Z^n$, and set 
\begin{eqnarray*}
&&
\Sigma^+ = \Sigma^+_s \cup \Sigma^+_\ell, \\
&&
\Sigma^+_s = \set{e_i - e_j, \; e_i + e_j}{1 \leq i < j \leq n}, \quad \Sigma^+_\ell = \set{2e_i}{1 \leq i \leq n},
\end{eqnarray*}
where $e_i$ is the $i$-th unit vector in $\Z^n, \; 1 \leq i \leq n$.
The set 
\begin{eqnarray*}
\Sigma_0 = \set{e_i - e_{i+1}}{1 \leq i \leq n-1} \cup \{2e_n \}
\end{eqnarray*}
forms the set of simple roots, and we denote by $\Delta$ the set of reflections associated with elements in $\Sigma_0$. Then 
$$
\Delta = \set{\tau_i = (i \; {i+1}) \in S_n}{1 \leq i \leq n-1} \cup \{ \tau \},
$$ 
which generates $W$.
For each $\sigma \in W$, we set
\begin{eqnarray} \label{Sigma-sigma} 
&&
\Sigma^+_s(\sigma) = \set{\alp \in \Sigma^+_s}{-\sigma(\alp) \in \Sigma^+}. 
\end{eqnarray}
We consider a pairing on $\Z^n \times \C^n$ given by
$$
\pair{t}{z} = \sum_{i=1}^n t_iz_i, \qquad (t \in \Z^n, \; z \in \C^n),
$$
which gives a $W$-invariant pairing on $\Sigma \times \C^n$, i.e., 
\begin{eqnarray*} \label{W-invariance}
\pair{\alp}{z} = \pair{\sigma(\alp)}{\sigma(z)}, \qquad (\alp \in \Sigma, \; z \in \C^n, \; \sigma \in W).
\end{eqnarray*}

\begin{thm} \label{th feq}
For $T \in \calH_n^{nd}$ and $\sigma \in W$, the spherical function $\omega_T(x;z)$ satisfies the following functional equation
\begin{eqnarray} \label{Gamma-sigma}
\omega_T(x; z) = \Gamma_\sigma(z) \cdot \omega_T(x; \sigma(z)),
\end{eqnarray}
where 
\begin{eqnarray*} \label{sigma-factors}
\Gamma_\sigma(z) = \dprod{\alp \in \Sigma^+_s(\sigma)}\,
\frac{1 - q^{\pair{\alp}{z}-1}}{q^{\pair{\alp}{z}} - q^{-1}},
\end{eqnarray*}
and we understand $\Gamma_\sigma(z) = 1$ if $\Sigma^+_s(\sigma) = \emptyset$.
In particular, the Gamma factor $\Gamma_\sigma(z)$ does not depend on $x$ nor $T$. 

\end{thm}

\slit
We note here that the factor $\pair{\alp}{z}$ for $\alp = e_i \pm e_j\; (i < j)$ in $s$-variable:
\begin{eqnarray} \label{alp-z}
\pair{\alp}{z} = 
\left\{ \begin{array}{ll}
  -(s_i+\cdots+s_{j-1})
     & \mbox{if } \alp = e_i-e_j\\[4mm]
 -(s_i+\cdots+s_{j-1}+ 2(s_j+\cdots + s_n))
     & \mbox{if } \alp = e_i+e_j 
\end{array}
 \right. .  
\end{eqnarray}

\slit
{\it Proof of Theorem~\ref{th feq}}. 
We define the Gamma factor $\Gamma_\sigma(z)$ by the equation (\ref{Gamma-sigma}). Then it is a rational function of $q^{z_1}, \ldots, q^{z_n}$ since $\omega_T(x; z)$ and $\omega_T(x; \sigma(z))$ are those functions, and Gamma factors satisfy the cocycle relations
\begin{eqnarray} \label{cocycle}
\Gamma_{\sigma_2\sigma_1}(z) = \Gamma_{\sigma_2}(\sigma_1(z)) \cdot \Gamma_{\sigma_1}(z), \qquad (\sigma_1, \sigma_2 \in W).
\end{eqnarray}
For convenience we set for $\alp \in \Sigma$
\begin{eqnarray} \label{f_alp}
f_{\alp}(\pair{\alp}{z}) = \left\{ \begin{array}{ll} 
1 & \mbox{ if } \alp = \pm 2e_i, \; (1 \leq i \leq n)\\
\displaystyle{\frac{1 - q^{\pair{\alp}{z}-1} }{ q^{\pair{\alp}{z}} -q^{-1} }}  & \mbox{ otherwise .} 
\end{array}
\right.
\end{eqnarray}
For an element $\sigma \in \Delta$ associated with some $\alp_0 \in \Sigma_0$, 
$$
\Sigma^+_s(\sigma) = \left\{ \begin{array}{ll}
\{ \alp_0 \} & \mbox{if } \alp_0 \in \Sigma^+_s\\
\emptyset & \mbox{if } \alp_0 \in \Sigma^+_\ell\; (\mbox{i.e., } \alp_0 = 2e_n),
\end{array} \right.
$$
and, by Remark~\ref{sect1 rem}, Remark~\ref{rem size1} and Theorem~\ref{th: tau}, 
\begin{eqnarray*} 
\Gamma_\sigma(z) =  f_{\alp_0}(\pair{\alp_0}{z}),
\end{eqnarray*}
which is independent of $x$ nor $T$. In general, assume that $\sigma \in W$ has the following shortest expression  
$$
\sigma = \sigma_{\ell} \cdots \sigma_{1}, 
$$
where $\sigma_i \in \Delta$ is the reflection associated with $\alp_i \in \Sigma_0$.
Then we see 
\begin{eqnarray*} \label{sigma-positive }
\set{\alp \in \Sigma^+ }{ \sigma(\alp) < 0} = 
\{\alp_1\} \cup \set{ \sigma_{1}\cdots \sigma_{k-1}(\alp_{k})   }{2 \leq k \leq \ell} .
\end{eqnarray*}
By using (\ref{cocycle}), (\ref{f_alp}) and the $W$-invariancy of the pairing $\pair{\;}{\;}$, we obtain
\begin{eqnarray*}
\Gamma_\sigma(z) &= &
\Gamma_{\sigma_{\ell}}(\sigma_{\ell-1} \cdots \sigma_1(z)) 
\cdots \Gamma_{\sigma_2}(\sigma_1(z)) \cdot \Gamma_{\sigma_1}(z)\\
&=&
f_{\alp_\ell}(\pair{\alp_\ell}{\sigma_{\ell-1} \cdots \sigma_1(z)}) \cdots
f_{\alp_2}(\pair{\alp_2}{\sigma_1(z)}) \cdot f_{\alp_1}(\pair{\alp_1}{z})\\
&=&
f_{\alp_\ell}(\pair{\sigma_1\cdots \sigma_{\ell-1}(\alp_{\ell})}{z}) \cdots
f_{\alp_2}(\pair{\sigma_1(\alp_2)}{z}) \cdot f_{\alp_1}(\pair{\alp_1}{z})\\
&=&
\prod_{\alp \in \Sigma^+_s(\sigma)}\, f_\alp(\pair{\alp}{z}),
\end{eqnarray*}
which completes the proof.
\qed

\bigskip
We will use the following explicit value $\Gamma_\rho(z)$ for a particular $\rho \in W$ in \S 5. 

\begin{cor} \label{cor rho}
Set $\rho \in W$ by
\begin{eqnarray*}
\rho(z_1, \ldots, z_n) = (-z_n, -z_{n-1}, \ldots, -z_1).
\end{eqnarray*}
Then
\begin{eqnarray*}
\Gamma_\rho(z) = 
\prod_{1 \leq i < j \leq n}\, 
\frac{1 - q^{z_i+z_j-1}}{q^{z_i+z_j}-q^{-1}}. 
\end{eqnarray*}
\end{cor}

\proof
Since 
$$
\Sigma^+_s(\rho) = \set{e_i + e_j}{1 \leq i < j \leq n},
$$
the assertion follows from Theorem~\ref{th feq}.
\qed

\begin{rem}{\rm
The above $\rho$ gives the functional equation of the hermitian Siegel series (cf.~\S 5), and it is interesting that such $\rho$ corresponds to the unique automorphism of the extended Dynkin diagram of the root system of type $(C_n)$, which was pointed out by Y.~Komori. 
}%
\end{rem}

\bigskip
\noindent
{\bf 2.4.} 
By using the functional equations (Theorem~\ref{th feq}) and the previous results on hermitian forms (Proposition~\ref{herm}), we obtain the following theorem, which gives us the information of the location of possible poles and zeros.

\begin{thm} \label{th: W-inv} 
Set 
\begin{eqnarray*}
G(z) = \prod_{\alp \in \Sigma^+_s}\, \frac{1 + q^{\pair{\alp}{z}}}{1 - q^{\pair{\alp}{z}-1}}.
\end{eqnarray*}
Then, for any $T \in \calH_n^{nd}$, the function $G(z) \cdot \omega_T(x; z)$ is holomorphic for all $z$ in $\C^n$ and $W$-invariant. In particular it is an element in $\C[q^{\pm z_1}, \ldots, q^{\pm z_n}]^{W}$.
\end{thm}

We denote by $\SKXT$ the subspace of $\CKXT$ consisting of compactly supported functions, which can be regarded as functions on $\frX_T$ of compactly supported functions modulo $U(T)$ on $\frX_T$ modulo.
Keeping the relation (\ref{change of var}) for $s$ and $z$, we consider the following integral
\begin{eqnarray} \label{Phi(xi,z)}
\Phi_T(z; \xi) = \int_{X_T^{op}}\, \xi(x) \abs{f_T(x)}^{s+\ve}dx, \qquad (\xi \in \SKXT)
\end{eqnarray}
where $dx$ is the $G$-invariant measure on $X_T$, and the right hand side is absolutely convergent for\begin{eqnarray*}
s \in \calD_0 
& = & 
\set{s \in \C^n}{ \real(s_i) \geq -\real(\ve_i), \; 1 \leq i \leq n}\\
& = &
\set{z \in \C^n}{-\frac12 \geq \real(z_n), \; \real(z_{i+1}) \geq \real(z_i) + 1, \; (1 \leq i \leq n-1)}.
\end{eqnarray*}
When $\xi$ is the characteristic function of $Kx$, $\Phi_T(z; \xi)$ is a constant multiple of $\omega_T(x;z)$, and any $\xi$ in $\SKXT$ is a finite linear sum of those characteristic functions. Thus we see that $\Phi_T(z; \xi)$ is a rational function of $q^{z_1}, \ldots, q^{z_n}$ and satisfy the same functional equations for $\omega_T(x; z)$, i.e.,
\begin{eqnarray} \label{feq of Phi}
\Phi_T(z; \xi) = \Gamma_\sigma(z) \cdot \Phi_T(\sigma(z); \xi), \qquad (\sigma \in W, \; \xi \in \SKXT).
\end{eqnarray}
Since $G(\sigma(z)) = G(z) \cdot \Gamma_\sigma(z)$ for $\sigma \in \Delta$, we see $G(z) \cdot \Phi_T(z; \xi)$ is invariant under the action of $\Delta$, hence it is $W$-invariant by cocycle relations. 
Since $G(z)$ is holomorphic for $z \in \calD_0$, 
we see $G(z)\cdot \Phi_T(z; \xi)$ is holomorphic for 
$$
z \in \bigcup_{\sigma \in W}\, \sigma(\calD_0).
$$
On the other hand, in a similar manner to the proof of Theorem~\ref{aim S_n}, we see
\begin{eqnarray*}
\Phi_T(z; \xi) = \int_{X_T^{op}}\, \xi(x) \zeta^{(n)}(D(x); s) dx,
\end{eqnarray*}
where $D(x) = x_2\cdot T^{-1}$ and $\zeta^{(n)}(y; s)$ is the spherical function on $\calH_n^{nd}$ (cf. the proof of Theorem~1.5), and recall that $G_1(z) \cdot \zeta^{(n)}(y; z)$ is holomorphic for $z \in \C^n$. 
Setting
\begin{eqnarray*}
G(z) = G_1(z)\cdot G_2(z), \qquad G_2(z) = \prod_{1 \leq i < j \leq n}\, \frac{1 + q^{z_i+z_j}}{1- q^{z_i+z_j-1}},
\end{eqnarray*}
we see $G(z) \cdot \Phi_T(z; \xi)$ is holomorphic for 
$$
z \in \calD_1 = \set{z \in \C^n}{\real(z_i+z_j) \ne 1, \; 1 \leq i < j \leq n},
$$
since $G_2(z)$ is holomorphic for $z \in \calD_1$ and $\xi$ is compactly supported.
Since $G(z)\cdot \Phi_T(z; \xi)$ is $W$-invariant, it is holomorphic for 
$$
z \in \wt{\calD} = \bigcup_{\sigma \in W}\, \sigma(\calD_0 \cup \calD_1).
$$
Since $\wt{\calD}$ is connected, $G(z)\cdot \Phi_T(z; \xi)$ is holomorphic in the convex hull $\C^n$ of $\wt{\calD}$.

Taking the characteristic function of $Kx$ for $\xi$, we obtain the theorem. 
\qed

\vspace{2cm}
\newcommand{\mapdownr}[1]{\Big\downarrow
   \rlap{$\vcenter{\hbox{$\scriptstyle#1$}}$ }}
\newcommand{\gen}[1]{\langle#1\rangle}

\setcounter{section}{2}
\setcounter{equation}{0}

\section{Explicit formulas}   
{\bf 3.1.}
Set
\begin{eqnarray}
\Lam_n^+ = \set{\lam \in \Z^n}{\lam_1 \geq \lam_2 \geq \cdots \geq \lam_n \geq 0},
\end{eqnarray}
and, for each $\lam \in \Lam_n^+$, 
\begin{eqnarray}
&&
\pi^\lam = Diag(\pi^{\lam_1}, \ldots, \pi^{\lam_n}) \in \calH_n^{nd}, \nonumber \\
&&
x_\lam = \twovector{\xi \pi^\lam}{1_n} \in \frX_{\pi^\lam}, \nonumber\\
&&
\omega_\lam(x; z) = \omega_T(x; z) \quad \mbox{for }\; T = \pi^\lam.
\end{eqnarray}

Then we obtain

\begin{thm} \label{th: explicit}
For $\lam \in \Lam_n^+$, one has the following explicit expression:
\begin{eqnarray} 
\lefteqn{\omega_\lam(x_\lam; z)} \nonumber \\
 &=&
\frac{(-1)^{\sum_i \lam_i(n-i+1)}  q^{-\sum_i \lam_i(n-i+\frac12)} (1 - q^{-2})^n}{\prod_{i=1}^{2n}(1 - (-q^{-1})^i)}  
\times \frac{1}{G(z)} \times 
\sum_{\sigma \in W}\, q^{-\pair{\lam}{\sigma(z)}} H(\sigma(z)) ,
\nonumber
\end{eqnarray}
where $G(z)$ is the same as in Theorem~\ref{th: W-inv} and
\begin{eqnarray*}
H(z) &=& 
\prod_{\alp \in \Sigma_s^+}\, 
\frac{1 + q^{\pair{\alp}{z}-1}}{1 - q^{\pair{\alp}{z}}} 
\prod_{\alp \in \Sigma_\ell^+}\, 
\frac{1 - q^{\pair{\alp}{z}-1}}{1-q^{\pair{\alp}{z}}}. 
\end{eqnarray*}
\end{thm}

\begin{rem} \label{rem: explicit}
{\rm 
By Theorem~\ref{th: W-inv}, the main part 
$$
H_\lam(z) = \sum_{\sigma \in W}\, \sigma\left( q^{-\pair{\lam}{z}} H(z) \right) = \sum_{\sigma \in W}\,  q^{-\pair{\lam}{\sigma(z)}} H(\sigma(z)) 
$$ 
of $\omega_\lam(x_\lam; z)$ belongs to $\C[q^{\pm z_1}, \ldots, q^{\pm z_n}]^W$. 
 Further we see directly in a standard way that the set $\set{H_\lam(z)}{\lam \in \Lam_n^+}$ forms its $\C$-basis. 
On the other hand, $H_\lam(z)$ is a special case of $P_\lam$ (up to a scalar factor) introduced by I.~G.~Macdonald (\cite{Mac} \S 10) in a generous context of orthogonal polynomials associated with root systems. 
}
\end{rem}

We will prove the above theorem by using a general expression formula given in \cite{French} (or in \cite{JMSJ} ) of spherical functions on homogeneous spaces, which is based on functional equations of finer spherical functions and some data depending only on the group $G$. 
We need to check the assumptions there. Let $\G$ be a connected reductive linear algebraic group and $\X$ be an affine algebraic variety which is $\G$-homogeneous, where everything is assumed to be defined over a $p$-adic field $k$. For an algebraic set, we use the same ordinary letter to indicate the set of $k$-rational points. Let $K$ be a special good maximal compact open subgroup of $G$, and $\B$ a minimal parabolic subgroup of $\G$ defined over $k$ satisfying $G = KB = BK$. 
We denote by $\frX(\B)$ the group of rational character of $\B$ defined over $k$ and by $\frX_0(\B)$ the subgroup consisting of those characters associated with some relative $\B$-invariant on $\X$ defined over $k$. In these situation, the assumptions are the following:

\mslit
$(A1)$ $\X$ has only a finite number of $\B$-orbits (, hence there is only one open orbit).

\mslit
$(A2)$ A basic set of relative $\B$-invariants on $\X$ defined over $k$ can be taken by regular functions on $\X$.

\mslit
$(A3)$ For $y \in \X$ not contained in the open orbit, there exists some $\psi$ in  $\frX_0(\B)$ whose restriction to the identity component of the stabilizer $\bH_y$ of $\G$ at $y$ is not trivial.

\mslit
$(A4)$ The rank of $\frX_0(\B)$ coincides with that of $\frX(\B)$.

\slit
In the present situation, as is noted in Remark~\ref{rem: realization}, we may understand $\G = U(H_n)$ as an algebraic group defined over $k$, $G = \G(k)$, $B = \B(k)$ for the Borel subgroup defined over $k$, $K = \G(\calO_k)$, and $X = X_T$ as the set of $k$-rational points of the affine algebraic variety $\X = \frX_T/U(T)$, and we recall that relative invariants $f_{T,i}(x)$ and the spherical function $\omega_T(x; s)$ can be regarded as functions on $X_T$.

It is easy to see the present $(\X, \B)$ satisfies the conditions $(A1)$, $(A2)$ and $(A4)$ (cf. Lemma~1.1, (\ref{rel inv}) and (\ref{asso char}) ), in particular, the unique Zariski open $\B$-orbit is given by $\X^{op} = \set{x \in \X}{f_{T,i}(x) \ne 0, \; 1 \leq i \leq n}$ (cf. (\ref{Xop})). 
We admit the condition $(A3)$, which is proved in \S 3.2, and give a proof of Theorem~3.1. 

\mslit
We recall the notation in the proof of of Proposition~\ref{many sph}.  
By the functional equation of $\omega_T(x;z)$ (Theorem~\ref{th feq}), we have for each $\sigma \in W$
\begin{eqnarray} 
\omega_T(x; z_\chi) &=& 
\Gamma_\sigma(z_\chi) \omega_T(x; \sigma(z_\chi))\nonumber \\
&=&
\Gamma_\sigma(z_\chi) \omega_T(x; \sigma(z)_{\sigma(\chi)}), \label{feq for omega-chi}
\end{eqnarray}
by taking a suitable character $\sigma(\chi)$ of $\calU$. 
When $\chi$ is the trivial character ${\bf 1}$, the equation (\ref{feq for omega-chi}) coincides with the original functional equation of $\omega_T(x; z)$ and $\Gamma_\sigma(z_{\bf 1}) = \Gamma_\sigma(z)$.
By (\ref{omega(z_chi)}) and (\ref{feq for omega-chi}), we obtain vector-wise functional equations for finer spherical functions $\omega_{T,u}(x;z)$ 
\begin{eqnarray} \label{matrix-feq}
\left( \omega_{T,u}(x; z) \right)_{u \in \calU} = A^{-1} \cdot G(\sigma,z)\cdot \sigma A \cdot \left( \omega_{T,u}(x; \sigma(z)) \right)_{u \in \calU}, \qquad \sigma \in W,
\end{eqnarray}
where
\begin{eqnarray*}
A = (\chi(u))_{\chi, u}, \quad
\sigma A = (\sigma(\chi)(u))_{\chi, u} \in GL_{2^{n-1}}(\Z),
\end{eqnarray*}
$\chi$ runs over characters of $\calU$, $u \in \calU$, and $G(\sigma, z)$ is the diagonal matrix of size $2^{n-1}$ whose $(\chi, \chi)$-component is $\Gamma_\sigma(z_\chi)$. 
We denote by $U$ the Iwahori subgroup of $K$ compatible with $B$ and take the normalized Haar measure $du$ on $U$. 
It is easy to see 
\begin{eqnarray*}
Ux_\lam \subset B x_\lam \quad \mbox{and} \quad \abs{f_T(ux_\lam)}^s = \abs{f_T(x_\lam)}^s,
\end{eqnarray*}
which means $x_\lam \in \calR^+$ in the sense of (2.8) in \cite{French}.
We set
\begin{eqnarray} 
\delta_{u}(x_\lam, z) &=& \dint{U}\, \abs{f_T(u x_\lam)}_u^{s+{\ve}} du  = 
 \left\{ \begin{array}{ll} 
\abs{f_T(x_\lam)}^{s+{\ve}} & \mbox{if } x_\lam  \in \frX_{T,u}\\
{} & {}\\
0 & \mbox{otherwise}.
\end{array}
\right\} \nonumber \\[2mm]
&=&
\label{delta-u}
 \left\{ \begin{array}{ll} 
(-1)^{\sum_i \lam_i(n-i+1)}  q^{-\sum_i \lam_i(n-i+\frac12)} q^{-<\lam, z>} & 
\mbox{if } x_\lam  \in \frX_{T,u}\\
{} & {}\\
0 & \mbox{otherwise}.
\end{array}
\right.
\end{eqnarray}
Applying Theorem~2.6 in \cite{French} to our present case, we obtain
\begin{eqnarray} \label{gen-formula}
\left( \omega_{T,u}(x_\lam; z) \right)_{u \in \calU} = 
\frac{1}{Q} \sum_{\sigma \in W}\, \gamma(\sigma(z)) \left( A^{-1} \cdot G(\sigma, z) \cdot \sigma A \right) \left( \delta_u(x_\lam, \sigma(z)) \right)_{u \in \calU},
\end{eqnarray}
where
\begin{eqnarray*}
&&
Q = \sum_{\sigma \in W}\, [U\sigma U : U]^{-1} = 
\prod_{i = 1}^{2n}\, \left( 1-(-1)^iq^{-i} \right) \Big{/} (1 - q^{-2})^{n},\\
&&
\gamma(z) =  \prod_{\alp \in \Sigma_s^+}\, \frac{1 - q^{2\pair{\alp}{z}-2}}{1 - q^{2\pair{\alp}{z}}}
\cdot
\prod_{\alp \in \Sigma^+_\ell}\, \frac{ 1 - q^{\pair{\alp}{z}-1} }{ 1 - q^{\pair{\alp}{z}}}. 
\end{eqnarray*}
By (\ref{matrix-feq}), (\ref{delta-u}), (\ref{gen-formula}), and the orthogonal relation of characters, we obtain
\begin{eqnarray*}
\omega_T({x_\lam}; z) &=& \sum_{u \in \calU}\, {\bf 1}(u) \omega_{T,u}(x_\lam; z)\\
&=&
\frac{(-1)^{-\sum_i \lam_i(n-i+1)}  q^{-\sum_i \lam_i(n-i+\frac12)}}{Q} \times \sum_{\sigma \in W}\, {\gamma(\sigma(z))} \Gamma_\sigma(z)  q^{-<\lam, \sigma(z)>}.
\end{eqnarray*}
Since we have
\begin{eqnarray*}
&&
\Gamma_\sigma(z) = \frac{G(\sigma(z))}{G(z)} \quad \mbox{(by Theorem~\ref{th: W-inv})},\\
&&
\gamma(z) \cdot G(z) = \prod_{\alp \in \Sigma^+_s} \frac{1 + q^{\pair{\alp}{z}-1}}{1 - q^{\pair{\alp}{z}}} \times \prod_{\alp \in \Sigma^+_\ell} \frac{1 - q^{\pair{\alp}{z}-1}}{1 - q^{\pair{\alp}{z}}} = H(z),
\end{eqnarray*}
we obtain
\begin{eqnarray*}
\omega_T({x_\lam}; z) &=& 
\frac{(-1)^{\sum_i \lam_i(n-i+1)}  q^{-\sum_i \lam_i(n-i+\frac12)} )(1-q^{-2})^n}{\prod_{i=1}^{2n}(1 - (-q^{-1})^i)}  
\times \frac{1}{G(z)} 
\times 
\sum_{\sigma \in W}\, \sigma\left(q^{-\pair{\lam}{z}} H(z) \right),
\end{eqnarray*}
which proves the theorem.
\qed

\slit
By Theorem~\ref{th: explicit} and Proposition~\ref{T[h]}, we get the explicit formula of $\omega_T(x;s)$ at many points.
For $\lam \in \Lam_n^+$ and $T \in \calH_n^{nd}$, it is known that $T$ and $\pi^\lam$ belong to the same $GL_n(k')$-orbit in $\calH_n^{nd}$ 
if and only if 
$$
v_\pi(\det T) \equiv \abs{\lam} \pmod{2},
$$
where $\abs{\lam} = \sum_{i=1}^n\, \lam_i$. 
And then, there exists some $h_\lam \in GL_n(k')$ for which $\pi^\lam[h_\lam] = T$ and $x_\lam h_\lam \in \frX_T$. Hence we have the following.

\begin{thm} \label{th: many explicit}
Let $T \in \calH_n^{nd}$ and $\lam \in \Lam_n^+$ and assume that $v_\pi(\det T) \equiv \abs{\lam} \pmod{2}$. Taking $h_\lam \in GL_n(k')$ for which $\pi^\lam[h_\lam] = T$, one has $x_\lam h_\lam \in \frX_T$ and 
\begin{eqnarray*}
\omega_T(x_\lam h_\lam; z) &=&
\omega_\lam(x_\lam ; z) \\
&=&
\frac{(-1)^{\sum_i \lam_i(n-i+1)}  q^{-\sum_i \lam_i(n-i+\frac12)} (1 - q^{-2})^n}{\prod_{i=1}^{2n}(1 - (-q^{-1})^i)}  
\cdot \frac{1}{G(z)} \cdot 
\sum_{\sigma \in W}\, \sigma\left(q^{-\pair{\lam}{z}} H(z)\right) .
\end{eqnarray*}
Further, each of such $\lam$'s gives a different $K$-orbit
$$
Kx_\lam h_\lam U(T) \quad \mbox{in }\; K \backslash X_T \; \Big( = K \backslash \frX_T/U(T)  \Big).
$$
\end{thm}

The latter statement follows from the explicit formula, since different $\lam$ gives the different value $\omega_T(x_\lam h_\lam; z)$ as a rational function of $q^{z_1}, \ldots, q^{z_n}$.

\slit
\noindent
{\bf 3.2.} In this subsection we prove the present $(\X, \B)$ satisfies the condition $(A3)$. We consider the action of $G \times U(T)$ on $\frX_T$ defined by $(g,h) \circ x = gxh^{-1}$. Then, the stabilizer $B_y$ of $B$ at $yU(T) \in X_T$ coincides with the image $B_{(y)}$ of the projection to $B$ of the stabilizer $(B \times U(T))_y$ at $y \in \frX_T$ to $B$. Hence, in our case, the condition $(A3)$ is equivalent to the following:

\medskip
\noindent
$(C)$ : For each $y \in \frX_T$ not contained in $\frX_T^{op}$, there exists $\psi \in \frX(\B)$ whose restriction to the identity component of $B_{(y)}$ is not trivial.

\medskip
\noindent
It suffices to prove the condition $(A3)$ (or $(C)$) over the algebraic closure $\ol{k}$ of $k$,
since, for a connected linear algebraic group $\bH$, $\bH(k)$ is dense in $\bH(\ol{k})$.
Then, we need to consider only for the case $T = 1_n$, since $\frX_T$ is isomorphic to $\frX_{T[g]}$ by $x \longmapsto xg$ and $B_{(x)} = B_{(xg)}$ for $g \in GL_n$;  and for simplicity of notation, we write $f_i(x)$ instead of $f_{T,i}(x)$. 
Until the end of this subsection, we consider algebraic sets over $\ol{k}$, extend the involution $*$ on $k'$ to $\ol{k}$, indicate it by ${}^{\overline{\hspace*{2.3mm}}}$, and write $\ol{x} = (\ol{x_{ij}}) \in M_{\ell m}(\overline{k})$ for $x = (x_{ij}) \in M_{\ell m}(\overline{k})$. 

Then, our situation is the following: 
\begin{eqnarray*}
&&
\frX = \frX_{1_n} = \set{x \in M_{2n,n}}{H_n[x] = 1_n},\\
&&
\left( U(H_n) \times U(1_n) \right) \times \frX \longrightarrow \frX, \quad ((g,h), x) \longmapsto (g,h) \circ x = gxh^{-1}, \nonumber 
\end{eqnarray*}
and $B$ is the Borel subgroup of $U(H_n)$ (as in (1.3)). 
We introduce a $(GL_{2n}\times GL_n)$-set $\wt{\frX}$ as follows:
\begin{eqnarray} \label{wt(X)}
&& \label{bigger action}
\wt{\frX} = \set{(x,y) \in M_{2n,n} \oplus M_{2n,n}}{{}^ty H_n x = 1_n}  \\[2mm]
&&
(g,h) \star (x,y) = (gxh^{-1}, \dot{g}y{}^th), 
\qquad ((g,h) \in GL_{2n}\times G_n, \; \dot{g} = H_n{}^tg^{-1}H_n). \nonumber
\end{eqnarray}
We write an element of $\wt{\frX}$ as $(x, y) = (\twovector{x_1}{x_2}, \twovector{y_1}{y_2})$ with $x_i, y_i \in M_n$, then the above condition is the same with
$$
{}^tx_1 y_2 + {}^tx_2 y_1 = 1_n. 
$$
We fix a Borel subgroup $P$ of $GL_{2n}$ by 
$$
P = \set{\twomatrix{p}{r}{0}{q} \in GL_{2n}}{{}^tp, \, q \in B_n, \; \; r \in M_n},
$$
where $B_n$ is the Borel subgroup of $GL_n$ consisting of the lower triangular matrices.
The involution
$g \longmapsto \dot{g} = H_n{}^tg^{-1}H_n$ on $GL_{2n}$  
induces an involution on $P$ : 
\begin{eqnarray}  \label{dot-p}
\twomatrix{p}{r}{0}{q} \longmapsto \twomatrix{{}^tq^{-1}}{-{}^tq^{-1}\, {}^tr\, {}^tp^{-1}}{0}{{}^tp^{-1}}.
\end{eqnarray}

Since $\dot{g} = \overline{g}$ for $g \in U(H_n)$ and ${}^th = \overline{h}^{-1}$ for $h \in H(1_n)$, the embedding $\iota : \frX \longmapsto \wt{\frX}, \; x \longmapsto (x, \overline{x})$ is compatible with the actions, i.e., we have the commutative diagram 
\begin{equation*}
\begin{array}{ccccc}
\left( U(H_n) \times U(1_n) \right) &\times& \frX & \stackrel{\circ}{\longrightarrow} & \frX\\
\mapdownr{incl.} & {} & \mapdownr{\iota} & {} & \mapdownr{\iota}\\[2mm]
\left( GL_{2n} \times GL_n \right) &\times& \wt{\frX} & \stackrel{\star}{\longrightarrow} & \wt{\frX}.
\end{array}
\end{equation*}
For $(x,y) \in \wt{\frX}$ and $p \in P$, set
\begin{eqnarray} \label{rela-inv wt(X)}
\wt{f_i}(x,y) = d_i(x_2{}^ty_2), \quad \wt{\psi}_i(p) = \prod_{1 \leq j \leq i}\, p_j^{-1}p_{n+j}, \quad (1 \leq i \leq n),
\end{eqnarray}
where 
$p_j$ is the $j$-th diagonal component of $p$.
Then $\wt{f}_i(x,y)$'s are relative $P$-invariants on $\wt{\frX}$ associated with characters $\wt{\psi}_i$, $\wt{f}_i(x,\overline{x}) = f_i(x)$ for $x \in \frX$, and $\wt{\psi}_i \vert_B = \psi_i$. We set 
$$
\calS = \set{(x,y) \in \wt{\frX} \cap (P \times GL_n)\star \iota(\frX)}{\prod_{i=1}^n\, \wt{f}_i(x,y) = 0}.
$$

For $\alp = (x, y) \in \wt{\frX}$, we denote by $H_\alp$ the stabilizer of $P \times GL_n$ at $\alp$, and by $P_\alp$ the identity component of the image of $H_\alp$ by the projection to $P$.
In order to prove the condition $(C)$, it suffices to show the following:

\medskip
\noindent
$(\wt{C}):$ For each $\alp \in \calS$,  there exists some $\psi \in \gen{\wt{\psi}_i \mid 1 \leq i \leq n}$ whose restriction to $P_\alp$ is not trivial.

\medskip
We have only to consider $(\wt{C})$ for representatives under the action of $P \times GL_n$. 
In the following we consider the case $n \geq 2$, since $\frX_T = \frX_T^{op}$ for $n = 1$ and there is nothing to prove.
We denote by $\delta_i(a) \in GL_n$ the diagonal matrix whose $j$-th entry is $1$ except the $i$-th which is $a \in GL_1$.

\begin{lem}
The condition $(\wt{C})$ is satisfied for $(x,y) \in \calS$ for which $\det x_2 \ne 0$ or $\det y_2 \ne 0$.
\end{lem}

\proof
Let $\alp = (x,y) \in \calS$ and $\det x_2 \ne 0$. 
Then by the action of $P \times GL_n$, we may assume that $x_2 = 1_n $ and $x_1 = 0$, then $y_1 = 1_n$ since ${}^tx H_n y = 1_n$. Since $\alp \in (P \times GL_n)\star \iota(\frX)$, $y_2$ can be written as 
$$
y_2 = p h, \qquad (p \in B_n, \; h \in GL_n, \, {}^t\overline{h}= h),
$$
and $0 = \prod_i\wt{f}_i(\alp) = \prod_i d_i(y_2) = \prod_i d_i(h)$.
For $q \in B_n$, we have
\begin{eqnarray*}
(\twomatrix{{}^t\,q^{-1}{}^tp}{0}{0}{q}, \, q) \star \alp &=& 
(\twomatrix{{}^tq^{-1}\, {}^tp}{0}{0}{q}\twovector{0}{1_n} q^{-1}, \;
\twomatrix{{}^tq^{-1}}{0}{0}{q p^{-1}}\twovector{1_n}{y_2}{}^tq)\\
&=&
( \twovector{0}{1_n}, \twovector{1_n}{qh{}^tq}) \left( = \beta, \; \mbox{say} \right).
\end{eqnarray*}
Hence, by taking a suitable $q \in B_n$, we may make $qh{}^tq = 1_r \bot h_1, \; 0 \leq r < n$ such that $h_1$ is a hermitian matrix satisfying\\

\quad \begin{minipage}{14.5cm}
$\cdot$ the first row and column are zero, \quad or

$\cdot$ for some $i, \; (1 < i \leq  n-r)$, each entry in the first row and column or in the $i$-th row and column is $0$ except at $(1,i)$ or $(i,1)$ which are $1$. 
\end{minipage}\\

\noindent
Then $H_{\beta}$ contains the following elements, according to the above type of $h_1$, 
\begin{eqnarray*}
(\twomatrixplus{\delta_{r+1}(a)}{}{}{1_n}, 1_n) \quad \mbox{or}\quad
(\twomatrixplus{\delta_{r+1}(a)}{}{}{\delta_{r+i}(a)}, \delta_{r+i}(a))\quad (a \in GL_1),
\end{eqnarray*}
and we see $\wt{\psi}_{r+1} \not\equiv 1$ on $P_{\beta}$.

The case $\alp = (x, y) \in \calS$ with $\det y_2 \ne 0$ is reduced to the case $\det x_2 \ne 0$, since $\beta = (y,x) \in \calS$, $H_\beta = \set{(\dot{p}, {}^th^{-1})}{(p,h) \in H_\alp}$ and $\wt{\psi}_i(\dot
p) = \wt{\psi}_i(p)^{-1}$.
\qed

\bigskip
Now we have to consider for $(x,y) \in \calS$ such that $\det x_2 = \det y_2 = 0$.
We set
\begin{eqnarray*}
&&
\calS_0 = \set{(x,y) \in \calS}{\det x_2 = \det y_2 = 0},\\
&&
J(i_1, i_2, \ldots, i_t) \in M_{n t}\,; \;  
\begin{array}{l} 1 \leq i_1 < i_2 < \cdots < i_t \leq n,  \\
\mbox{the entry at $(i_j, j)$ is $1$, and all the other entries are 0.}
\end{array}
\end{eqnarray*}

\begin{lem} \label{shape}
By the action of $P \times GL_n$, every element in $\calS_0$ becomes the following type, 
\begin{eqnarray*} \label{basic form}
&&
(\twomatrixplus{0}{J_1}{J_2}{0}, \twomatrixplus{z_1}{0}{z_2}{z_3} ),  \qquad (J_1, z_3 \in M_{n \ell}, \; J_2, z_1, z_2 \in M_{n k}),
\end{eqnarray*}
where
\begin{eqnarray*}
J_1 = J(r_1, r_2,\ldots,r_\ell), \quad J_2 = J(e_1, e_2,\ldots,e_k), \quad 1 \leq \ell, k < n, \; \; \ell + k = n, 
\end{eqnarray*}
and 
\begin{eqnarray*} 
&&
\mbox{the $e_j$-th row of $z_1$ is the same as in $J_2$ and $(i, j)$-entry is $0$ if $i < e_j, \; 1 \leq j \leq k$,} \qquad \nonumber\\
&&
\mbox{the $r_j$-th row of $z_2$ is $0, \; 1 \leq j \leq \ell$,} \label{assump-z} \\
&&
\mbox{the $r_j$-th row of $z_3$ is the same as in $J_1$ and $(i, j)$-entry is $0$ if $i > r_j,\; 1 \leq j \leq \ell$.}\qquad \nonumber
\end{eqnarray*}
\end{lem}

\proof
Take an $\alp = (x, y) \in \calS_0$ and let ${\rm rank}(x_2) = k$. Then $1 \leq k <n$, and by the action of $P \times GL_n$, we make $x$ into  
$$
\twomatrixplus{0}{x'}{J_2}{0}.
$$
Then, the rank of $x'$ must be $\ell = n-k$, since $x \in \wt{\frX}$, and we may make $x'$ into $J_1$, i.e. $x$ into the required type.
Further, the $e_j$-th rows in $y_1$ must be the same as in $(J_2 \mid 0)$ and the $r_j$-th rows in $y_2$ must be the same as in $(0 \mid J_1)$. 

Multiplying $y$ by a suitable element $p \in P$ from the left we may make the latter $\ell$ columns of $y_1$ into $0$ and $(i, k + j)$-entry of $y_2$ for $1 \leq j \leq \ell, \; i > r_j$ into $0$, while $\dot{p}x = x$. 
Since $(e_j, r)$-entry of $y_1$ is $0$ unless $r = j$, we may make $(i, j)$-entry of $y_1$ for $1 \leq j \leq k, \; i < r_j$ into $0$ as keeping $x$. Thus we obtain a matrix of the form as in the statement.
\qed

\begin{lem}
The condition $(\wt{C})$ is satisfied for elements in $\calS_0$.
\end{lem}

\proof 
We may assume $\alp = (x, y) \in \calS_0$ has the form as in Lemma~\ref{shape}. Then, for any $a \in GL_1$,  
\begin{eqnarray*}
\begin{array}{ll}
(\twomatrixplus{1_n}{0}{0}{\delta_1(a)}, 1_n) \in H_\alp & \mbox{ if }\; e_1 > 1, \\[3mm]
(\twomatrixplus{\delta_1(a)}{0}{0}{1_n}, \delta_{k+1}(a)) \in H_\alp & \mbox{ if }\; r_1 = 1,\\[3mm]
(\twomatrixplus{a1_n}{0}{0}{1_n}, a1_n) \in H_\alp & \mbox{ if }\; z_2 = 0.
\end{array}
\end{eqnarray*}
When $e_1 = 1$, $r_1 > 1$ and $z_2 \ne 0$, we modify $\alp$ into $\beta = (x, y')$ by the $P \times GL_n$-action as below:
\begin{eqnarray} 
&&
y ' = \twomatrixplus{z'_1}{0}{z'_2}{z'_3},\; z'_2 \ne 0, \nonumber \\
&&
\mbox{the $r_j$-th row of $z'_3$ is the $j$-th unit vector (the same as in $J_1$) for  $1 \leq j \leq \ell$.} \nonumber \\
&&
\mbox{if the $i$-th row of $z'_2$ is not $0$, then the $i$-th row of $z'_3$ is $0$, $1 \leq i \leq n$}.
\label{finer form}
\end{eqnarray}
Then, for any $a \in GL_1$,
$$
(\twomatrixplus{D_n(a_i)}
{0}{0}{1_n}, \; \twomatrixplus{1_k}{0}{0}{a1_\ell}) \in H_\beta, 
$$
where
$D_n(a_i) = Diag(a_1, \ldots, a_n)$ with 
$$
a_i = \left\{ \begin{array}{ll} 
a & \mbox{if the $i$-th row of $z'_2$ is $0$}\\
1 & \mbox{if the $i$-th row of $z'_2$ is not $0$},
\end{array} \right. .
$$
Hence, for any $\alp \in \calS_0$, $\wt{\psi}_n \not\equiv 1$ on $P_\alp$.

Now we explain how to obtain $\beta$ as in (\ref{finer form}) from $\alp$ with $e_1=1, \; r_1 > 1$ and $z_2 \ne 0$. 
Let $k'$ be the rank of $x_2$. Then for suitable $p_0 \in B_n$, we make $z'_2 = p_0z_2$ such that
\begin{eqnarray}
&&
\mbox{there exist integers $1 \leq s_1 < s_2 < \cdots < s_{k'} \leq n$ such that }
\mbox{the $i$-th rows are $0$ except  } \nonumber \\
&&
\mbox{for } i \in \{s_1, \ldots, s_{k'} \}, \mbox{ and for each $i, \; 1 \leq i \leq k'$, there exists distinct $j_i$ for which } \nonumber \\
&&
\qquad \mbox{$1$ at $(s_i, j_i)$-entry}, \nonumber \\
&&
\qquad \mbox{$0$ at $(s_i, j)$-entry for $j < j_i$ and 
the $(i', j_i)$-entry for  $i' > s_i$. } \label{new z2}
\end{eqnarray}
Since every $r_i$-th row of $z_2$ is $0$, we may assume each $r_i$-th row of $p_0$ is the $r_i$-th unit vector, hence ${}^tp_0^{-1}J_1 = J_1$ and the $r_j$-th row of $p_0z_3$ remains to be the $j$-th unit vector. By a suitable matrix 
$$
h = \twomatrix{1_k}{C}{0}{1_\ell} \in GL_n,
$$
we make each $s_i$-th row of $z_3' = z_2'C + p_0z_3$ into $0$ for $1 \leq  i \leq k'$ and remain the other rows as the same as in $p_0z_3$. 
Take the matrix $D \in M_n$ by putting the $j$-th row of $z_1C$ into the $r_j$-th row for $1 \leq j \leq \ell$, and $0$ at all other entries, then $z_1C = Dz_3'$.  
Setting 
$$
p_1 = \twomatrix{1_n}{0}{0}{p_0}  \quad\mbox{and}\quad p = \twomatrix{1_n}{-D}{0}{1_n}, 
$$
we obtain 
$$
p p_1 yh = \twomatrix{1_n}{-D}{0}{1_n}\twomatrix{z_1'}{z_1C}{z_2'}{z_3'} = 
\twomatrix{z_1'}{0}{z_2'}{z_3'}, \qquad (z_1' = z_1-Dz_2').
$$
On the other hand, we have
\begin{eqnarray*}
\dot{p} \dot{p_1} x {}^th^{-1} &=& \twomatrix{1_n}{{}^tD}{0}{1_n} \twomatrix{{}^tp_0^{-1}}{0}{0}{1_n} \twomatrix{0}{J_1}{J_2}{0} 
\twomatrix{1_k}{0}{-{}^tC}{1_\ell} =   
\twomatrix{{}^t DJ_2 - J_1 {}^tC}{J_1}{J_2}{0},
\end{eqnarray*}
and ${}^tDJ_2$ and $J_1 {}^tC$ may have nonzero rows only at the $r_i$-th, $1 \leq i \leq \ell$, and 
\begin{eqnarray*}
(r_i, j)\mbox{-entry of } {}^tDJ_2 & =&
(e_j, i)\mbox{-entry of } z_1C = (j,i)\mbox{-entry of } C \\
&=& 
(r_i, j)\mbox{-entry of } J_1\, {}^tC.
\end{eqnarray*}
Thus we have the required element  
$$
\beta = (\dot{p}\dot{p_1}, {}^th) \star \alp = (\twomatrix{0}{J_1}{J_2}{0}, \twomatrix{z_1'}{0}{z_2'}{z_3'}).
$$
\qed

\bigskip
Thus we have shown the condition $(\wt{C})$ is satisfied for every $(x,y) \in \calS$, which shows that our $(\X, \B)$ satisfies the condition $(A3)$ and Theorem~3.1 is established.

\vspace{2cm}
\setcounter{section}{3}
\setcounter{equation}{0}

\section{Spherical Fourier transform on $\calS(K \backslash X_T)$} 
%
We consider the subspace $\SKXT$ of $\CKXU$ consisting of compactly supported modulo $U(T)$ functions, which is an $\hec$-submodule (cf. (\ref{conv prod})). 
We define the spherical Fourier transform $F_T$ on $\SKXT$, by setting
\begin{eqnarray}
&&
F_T : \SKXT \longrightarrow \C(q^{z_1}, \ldots, q^{z_n}), \nonumber \\ 
&&\qquad
\xi \longmapsto F_T(\xi)(z) 
= \int_{X_T} \xi(x) \Psi_T(x; z) dx,
\end{eqnarray}
where $\Psi_T(x;z) = G(z) \cdot \omega_T(x;z)$ and $dx$ is the $G$-invariant measure on $X_T$.
Since $\SKXT$ is spanned by the characteristic functions of double cosets $KxU(T)$ in $K\backslash \frX_T/U(T) = K \backslash X_T$, the image of $F_T$ is spanned by the set $\set{\Psi_T(x; z)}{ x \in \frX_T}$ over $\C$, and contained in 
\begin{eqnarray*}
\calR = \C[q^{\pm z_1}, \ldots, q^{\pm z_n}]^W
\end{eqnarray*}
by Theorem~\ref{th: W-inv}.
\newcommand{\bfe}{{\bf e}}
We decompose $\calR$ in the following 
\begin{eqnarray*}
&&
\calR = \bigoplus_{\bfe \in \{0, 1\}^n}\, s_1^{e_1}\cdots s_n^{e_n}\, \calR_0, 
\end{eqnarray*}
where
\begin{eqnarray*}
&&
\calR_0 = \C[q^{\pm 2z_1}, \ldots, q^{\pm 2z_n}]^W = \C[q^{2z_1}+q^{-2z_1}, \ldots, q^{2z_n}+q^{-2z_n}]^{S_n},
\end{eqnarray*}
and $s_i= s_i(z)$ is the $i$-th fundamental symmetric polynomial of $\set{q^{z_j}+q^{-z_j}}{1 \leq j \leq n}$; $\calR$ is a free $\calR_0$-module of rank $2^{n}$.
We set 
\begin{eqnarray*}
\calR_{even} = \bigoplus_{\bfe: even}\,  s_1^{e_1}\cdots s_n^{e_n}\, \calR_0, \quad
\calR_{odd} = \bigoplus_{\bfe: odd}\,  s_1^{e_1}\cdots s_n^{e_n}\, \calR_0,
\end{eqnarray*}
where $\bfe \in \{0, 1\}^n$ is even (resp. odd) if \, $\sum_{i=1}^n\, ie_i$ is even (resp. odd). For each $T \in \calH_n^{nd}$, we define
$$
\calR_{\langle{T}\rangle} 
$$
to be $\calR_{even}$ or $\calR_{odd}$ according to the parity of $v_\pi(\det(T))$.

\begin{thm} \label{sph tr}
For any $T \in \calH_n^{nd}$, one has a surjective $\hec$-module homomorphism
$$
F_T: \SKXT \longrightarrow \calR_{\langle{T}\rangle},
$$
and a commutative diagram
\begin{eqnarray}
\begin{array}{ccccc}
\hec & \times & \SKXT & \stackrel{*}{\longrightarrow} & \SKXT\\[2mm]
 {\Big\downarrow\llap{$\vcenter{\hbox{$\scriptstyle{\wr}\, $}}$ }}
& {} & {\Big\downarrow
   \rlap{$\vcenter{\hbox{$\scriptstyle{F_T}$}}$ }}
 & {\circlearrowleft} & 
 {\Big\downarrow\rlap{$\vcenter{\hbox{$\scriptstyle{F_T}$}}$ }}
\\[2mm]
\calR_0 &\times & \calR_{\langle{T}\rangle} &\longrightarrow& \calR_{\langle{T}\rangle},
\end{array}
\end{eqnarray}   
where the upper horizontal arrow is given by the action of $\hec$ on $\SKXT$, 
the left end vertical isomorphism is given by Satake isomorphism (\ref{satake iso})
$$
\hec \stackrel{\sim}{\longrightarrow} \calR_0, \; \phi \longmapsto \lam_z(\check{\phi}),  \qquad  (\check{\phi}(g) = \phi(g^{-1})),
$$
and the lower horizontal arrow is given by the ordinary multiplication in $\calR$. 
\end{thm}

\bigskip
\proof 
For $\phi \in \hec$ and $\xi \in \SKXT$, we have
\begin{eqnarray*}
F_T(\phi*\xi)(z) &=& \dint{X} \dint{G} \phi(g)\xi(g^{-1}x) dg \Psi_T(x; z) dx 
= \dint{X} \xi(y) \dint{G} \phi(g) \Psi_T(gy; z) dg dy \\[2mm]
&=& 
\dint{X} \xi(y) (\check{\phi} * \Psi_T(\;,z)) (y) dy = \lam_z(\check{\phi}) F_T(\xi)(z),  
\end{eqnarray*}
which gives the commutative diagram.

We recall the definition (\ref{def sph f}) of $\omega_T(x;z)$ and expand it in a region of absolute convergence. Then  
$$
\omega_T(x; z) = \sum_{\mu \in \Z^n}\, a_\mu q^{\pair{\mu}{z}},
$$
where $a_\mu = 0$ unless $\abs{\mu} \left( = \sum_{i=1}^n\, \mu_i \right) \equiv v_\pi(\det(T)) \pmod{2}$, since 
\begin{eqnarray*}
v_\pi(f_{T,n}(x)) &=& v_\pi(\det(x_2 T^{-1} x_2^*) \equiv v_\pi(\det(T)) \pmod{2}, \; \mbox{for any}\, x \in \frX_T^{op}\\ 
\pair{\mu}{z} &=& 
\sum_{i=1}^n\, \mu_iz_i = -\sum_{i=1}^n\, \mu_i(s_i+\cdots + s_n) \quad (\mbox{in $s$-variable})\\
&=&
-\mu_1s_1 - (\mu_1+\mu_2)s_2 - \cdots - (\mu_1+\cdots + \mu_n)s_n.
\end{eqnarray*}
Since
$$
G(z) = \prod_{i<j}\, \left( (1 + q^{z_i-z_j}+q^{z_i+z_i}+q^{2z_i}) \sum_{\ell, r \geq 0}\, q^{(\ell+r)z_i + (\ell - r)z_j - (\ell + r)} \right),
$$
can be expanded only in terms $q^{\pair{\nu}{z}}$ with $\abs{\nu}$ is even, we may expand   
$\Psi_T(x; z) = \omega_T(x; z) G(z)$ in terms $q^{\pair{\lam}{z}}$ with 
$\abs{\lam} \equiv v_\pi(\det(T)) \pmod{2}$,  hence 
%
\newcommand{\Image}{{\rm Im}}

\begin{eqnarray} \label{image}
\Image(F_T) \subset \calR_{\langle{T}\rangle}. 
\end{eqnarray}
%
%
%
On the other hand, by Remark~\ref{rem: explicit} and Theorem~\ref{th: many explicit} we see
$$
\Image(F_T) \supset \set{H_\lam(z)}{\lam \in \Lam_n^+, \; \abs{\lam}\equiv v_\pi(\det T)\pmod{2}},
$$
and the image of $F_T$ coincides with $\calR_{\langle{T}\rangle}$.
\qed

\bigskip

\begin{rem}{\rm 
We expect that the spherical Fourier transform $F_T$ is injective, which is equivalent to the identity 
\begin{eqnarray} \label{Cartan}
\frX_T = \bigcup_{\shita{\lam \in \Lam_n^+}{\abs{\lam} \equiv v_\pi(\det(T)) \pmod{2}}}\, K x_\lam h_\lam U(T),
\end{eqnarray}
where disjointness in the right hand side is known by Theorem~\ref{th: many explicit}.
If it is true, then $\SKXT$ would be a free $\hec$-module of rank $2^{n-1}$ and the set $\set{\Psi_T(x; z + \wt{u})}{ u \in \calU}$ would form a basis of spherical functions on $X_T$ corresponding to $z \in \C^n$ through $\lam_z$ (cf. Proposition~\ref{many sph}). 
This is true when $n = 1$ by Proposition~2.1, and we have the following. 
} 
\end{rem}

\bigskip

\begin{propos} 
Assume $n = 1$. Then the spherical transform $F_T$ is injective and $\SKXT$ is a free $\hec$-module of rank $1$, in fact the image coincides with 
$$
\C[q^{2z}+q^{-2z}] \; \mbox{if $v_\pi(T)$ is even}, \quad
(q^{z}+q^{-z})\C[q^{2z}+q^{-2z}] \; \mbox{if $v_\pi(T)$ is odd}.
$$ 
Any spherical function on $X_T$ corresponding to $z \in \C$ through $\lam_z$ is a constant multiple of $\omega_T(x;z)$.  
\end{propos}

\vspace{2cm}
\setcounter{section}{4}
\setcounter{equation}{0}
\section{An application to hermitian Siegel series}    
We recall the hermitian Siegel series, and give an integral representation and a new proof of the functional equation as an application of spherical functions. 

Let $\psi$ be an additive character of $k$ of conductor $\calO_k$.
For $T \in \calH_n(k')$ and $t \in \C$, the hermitian Siegel series $b_\pi(T;s)$ is defined by
\begin{eqnarray} \label{b-pi(t)}
b_\pi(T; t) = \dint{\calH_n(k')}\, \nu_\pi(R)^{-t} \psi(\tr(TR))dR,
\end{eqnarray}
where $\tr(\;)$ is the trace of matrix and $\nu_\pi(R)$ is defined as follows:
if the elementary divisors of $R$ with negative $\pi$-powers are $\pi^{-e_1}, \ldots, \pi^{-e_r}$,  then $\nu_\pi(R) = q^{e_1 + \cdots + e_r}$, and $\nu_\pi(R) = 1$ otherwise (cf. \cite{Shimura}-\S 13).
The right hand side of (\ref{b-pi(t)}) is absolutely convergent if $\real(t)$ is sufficiently large.

In the following we assume that $T$ is nondegenerate, since the properties of $b_\pi(T; t)$ can be reduced to the nondegenerate case. 
We give an integral expression of $b_\pi(T; t)$ in a similar argument for Siegel series in \cite{Arakawa}-\S 2.

We recall the set $\xx_T$ for $T \in \calH_n^{nd}(k')$
$$
\xx_T = \xx_T(k') = \set{x \in M_{2n,n}(k')}{H_n[x] = T}
$$
and take the measure $\abs{\Theta_T}$ on $\xx_T$ simultaneously as the fibre space of $T$ by the polynomial map $M_{2n,n}(k') \longrightarrow \calH_n(k'), x \longmapsto H_n[x]$ defined over $k$.
Then the following identity holds (cf. \cite{Ymzk}, \cite{Arakawa}-\S 2):
\begin{eqnarray*}
\lefteqn{\dint{\xx_T(k')}\, \phi(x) \abs{\Theta_T}(x)} \\
& =& \lim_{e \rightarrow \infty}\, \dint{\calH_n(\pi^{-e})}\, \psi(-\tr(Ty)) \dint{M_{2n,n}(k')}\, \phi(x) \psi(\tr(H_n[x]y)) dxdy, \nonumber
\end{eqnarray*}
where $\phi \in \calS(M_{2n,n}(k'))$, a locally constant compactly supported function on $M_{2n,n}(k')$, and $\calH_n(\pi^{-e}) = \calH_n(k') \cap M_n(\pi^{-e}\calO_{k'})$.

\bigskip
The following lemma can be proved in the similar line to the case of symmetric matrices (cf. \cite{Arakawa}-\S 2).

\begin{lem} \label{hss lem1}
If $\real(t)$ is sufficiently large, one has 
\begin{eqnarray} \label{lem hss}
\lefteqn{\dint{\xx_T(\calO_{k'})}\, \abs{N(\det x_2)}^{t-n} \abs{\Theta_T}(x)}\\
&=&
\lim_{e \longrightarrow \infty}\, \dint{\calH_n(\pi^{-e}\calO_{k'})}\, \psi(-\tr(Ty)) dy \dint{M_{2n,n}(\calO_{k'})}\, \abs{N(\det x_2)}^{t-n} \psi(\tr(H_n[x]y)) dx. \nonumber
\end{eqnarray}
\end{lem}

\bigskip
Denote by $\zeta(k'; t)$ the zeta function of the matrix algebra $M_n(k')$:
\[
\zeta(k'; t) = \int_{M_n(\calO_{k'})} \abs{\det x}_{k'}^{t-n}\,dx = \int_{M_n(\calO_{k'})} \abs{N(\det x)}^{t-n}\,dx,
\]
whose explicit formula is well-known:
\[
\zeta_n(k'; t) = \prod_{i=1}^n \frac{1-q^{-2i}}{1-q^{-2(t-i+1)}}. 
\]
Then we have the following integral expression of hermitian Siegel series.

\begin{thm} \label{hss thm2}
If $\real(t) > 2n$, we have
\begin{eqnarray*}
b_\pi(T; t) = \zeta_n(k';\frac{t}{2}) ^{-1} \times \dint{\xx_T(\calO_{k'})}\, 
\abs{N(\det x_2)}^{\frac{t}{2}-n} \abs{\Theta_T}(x).
\end{eqnarray*}
\end{thm}

\proof
We define the Fourier transform of $\phi \in \calS(M_n(k'))$ by 
$$
\what{\phi}(z) = \dint{M_n(k')}\, \phi(y) \psi(T_{k'/k}(\tr(yz^*)) dy,
$$
where $T_{k'/k}$ is the trace of the extension ${k'/k}$.
Since we have
\begin{eqnarray*}
\tr(H_n[x]y) = \tr(x_1^*x_2y) + \tr(x_2^*x_1y) = \tr(x_1^* (x_2y)) + \tr((x_2y)^* x_1) = T_{k'/k}(\tr(x_1(x_2y)^*),
\end{eqnarray*}
the second integral in the right hand side of (\ref{lem hss}) becomes
\begin{eqnarray*}
\lefteqn{\dint{M_n(\calO_{k'})}\, \abs{N(\det x_2)}^{t-n} \what{ch_{M_n(\calO_{k'})}}(x_2y)dx_2}\\
&=&
\dint{M_n(\calO_{k'})}\, \abs{N(\det x_2)}^{t-n} ch_{M_n(\calO_{k'})}(x_2y)dx_2\\
&=&
\dint{M_n(\calO_{k'})y^{-1} \cap M_n(\calO_{k'})}\, \abs{N(\det x_2)}^{t-n} dx_2\\&=&
\dint{M_n(\calO_{k'})D_y}\, \abs{\det x_2}_{k'}^{t-n} dx_2,
\end{eqnarray*}
where $D_y = 1_n$ if $y \in M_n(\calO_{k'})$, and $D_y = Diag(\pi^{e_1},\ldots,\pi^{e_r}, 1,\ldots,1)$ if the elementary divisors of $y$ with negative $\pi$-powers are $\pi^{-e_1},\ldots,\pi^{-e_r}$. 
Hence the second integral in the right hand side of (\ref{lem hss}) is equal to
\begin{eqnarray*}
\abs{\det D_y}_{k'}^{t} \dint{M_n(\calO_{k'})}\, \abs{\det x_2}_{k'}^{t-n} dx_2
= \nu_\pi(y)^{-2t} \times \zeta_n(k';t).
\end{eqnarray*}
Now by Lemma~\ref{hss lem1},  we obtain
\begin{eqnarray*}
\lefteqn{\dint{\xx_T(\calO_{k'})}\, \abs{N(\det x_2)}^{t-n} \abs{\Theta_T}(x)}\\
&=&
\zeta_n(k';t) \times \lim_{e \rightarrow \infty}\, \dint{\calH_n(\pi^{-e}\calO_{k'})}\,
\nu_\pi(y)^{-2t} \cdot \psi(-\tr(Ty))dy\\
& = &\zeta_n(k'; t) \times b_\pi(T; 2t),
\end{eqnarray*}
which gives the required identity.
\qed

\bigskip
Setting, in $s$-variable,
\begin{eqnarray} \label{s_t}
s_t = (1+\tfrac{\pi\sqrt{-1}}{\log q}, \ldots, 1+\tfrac{\pi\sqrt{-1}}{\log q}) + (0,\ldots, 0, \frac{t}{2} - n - \frac{1}{2}) \in \C^n,
\end{eqnarray}
we see 
\begin{eqnarray}
\dint{K}\, \abs{N(\det (kx)_2)}^{\frac{t}{2} - n}dk = 
\abs{\det T}^{\frac{t}{2} - n} \omega_T(x ; s_t).
\end{eqnarray}
Hence we may express $b_\pi(T;t)$ by using the spherical function $\omega_T(x; s)$.

\begin{propos} \label{prop2-1}
Denote the $K$-orbit decomposition of $\xx_T(\calO_{k'})$ as
$$
\xx_T(\calO_{k'}) = \sqcup_{i=1}^{r} Kx_i.
$$
Then one has
\begin{eqnarray*} \label{eq in prop}
b_\pi(T; t) &=& 
\abs{\det T}^{\frac{t}{2}-n}\, \prod_{i=0}^{n-1}\, (1 - q^{-t + 2i}) \times \sum_{i=1}^{r}\, c_i \cdot \omega_T(x_i; s_t), 
\end{eqnarray*}
where $c_i = \left( \prod_{i=1}^n ( 1 - q^{-2i}) \right)^{-1} \cdot vol(Kx_i)$.
\end{propos}

\proof
Since $\xx_T(\calO_{k'})$ is compact, it is a finite union of $K$-orbits, which we write as above.
By Theorem~\ref{hss thm2}, we have
\begin{eqnarray*}
\lefteqn{b_\pi(T; t) \times \zeta_n(k'; \frac{t}{2})}\\
&=&
\sum_{i=1}^r\, \dint{Kx_i}\, \abs{N(\det y_2)}^{\frac{t}{2}-n} \abs{\Theta_T}(y)\\
&=&
\sum_{i=1}^r\, \dint{Kx_i}\, \dint{K} \abs{N(\det (ky)_2)}^{\frac{t}{2}-n} dk \abs{\Theta_T}(y)\\
&=&
\abs{\det T}^{\frac{t}{2}-n}\, \sum_{i=1}^r\, c_i' \cdot \omega_T(x_i; s_t),
\end{eqnarray*}
where $c_i' = vol(Kx_i)$. Substituting the explicit value of $\zeta_n(k'; \frac{t}{2})$, we conclude the proof.
\qed

\slit
By using Theorem~\ref{th: W-inv}, we have the following.

\begin{cor} \label{cor}
The function $\{\prod_{i=0}^{n-1}\, (1 - (-1)^iq^{-t+i})\}^{-1} \times b_\pi(T; t)$ is holomorphic for any $t$, hence it is a polynomial in $q^{t}$ and $q^{-t}$.
\end{cor}

\proof
We denote by $z^*$ the corresponding value with $s_t$ in $z$-variable.
By Proposition~\ref{prop2-1} and Theorem~\ref{th: W-inv}, we see that
$$
b_\pi(T;t) = \prod_{i=0}^{n-1}\, (1 - q^{-t+2i})\cdot \frac{1}{G(z^*)} \times (\mbox{a holomorphic function}).
$$
By (\ref{s_t}), (\ref{alp-z}), and the definition of $G(z)$, we obtain
\begin{eqnarray*}
G(z^*) &\equiv&
\prod_{i<j}\, \frac{1 + (-1)^{j+i}q^{-t+i+j-1}}{1 - (-1)^{i+j}q^{-t + i + j -2}} \pmod{\C^\times}\\
&=& 
\prod_{i=1}^{n-1} \prod_{j=i+1}^n\, \frac{1 - (-1)^{i+j-1}q^{-t+i+j-1}}{ 1 - (-1)^{i+j-2}q^{-t+i+j-2}}\\
&=&
\prod_{i=1}^{n-1} \frac{1 - (-1)^{n+i-1}q^{n+i-1}}{1 + q^{-t + 2i-1}},
\end{eqnarray*}
and
\begin{eqnarray*}
\prod_{i=0}^{n-1}\, (1 - q^{-t+2i})\cdot \frac{1}{G(z^*)} \equiv \prod_{i=0}^{n-1}\, (1 - (-1)^iq^{-t+i}) \pmod{\C^\times},
\end{eqnarray*}
which completes the proof.
\qed

\begin{rem}  \label{shimura's decomp}
{\rm 
According to G.~Shimura \cite{Shimura} Theorem~13.6, one may express $b_\pi(T; t)$ as follows (including ramified hermitian and split cases):  
\begin{eqnarray} \label{shimura}
b_\pi(T; t) = f_T(q^{-t}) \cdot g_T(q^{-t}),
\end{eqnarray}
where $f_T(X)$ is an explicitly given  rational function of $X$, depending only on the type and size of $T$, and $g_T(X)$ is a (mysterious) polynomial with coefficients in $\Z$. For the unramified hermitian case, $f_T(X)$ is given for $T \in \calH_n^{nd}$ by
$$
f_T(X) = \prod_{i=0}^{n-1}(1- (-q)^iX), \qquad f_T(q^{-t}) = \prod_{i=0}^{n-1}(1- (-1)^iq^{-t+i}).
$$
In Corollary~\ref{cor}, we obtain the same factor $f_T(q^{-t})$ by using the spherical functions $\omega_T(x; z)$, and $f_T(q^{-t})^{-1} b_\pi(T; t)$ must be a polynomial in $q^{-t}$ with coefficients in $\Z$, which we don't see from $\omega_T(x; z)$.
} 
\end{rem}

\bigskip
Now we give the functional equation of the hermitian Siegel series by using the results of functional equations of the spherical functions $\omega_T(x; s)$.

\begin{thm} \label{th feq of Siegel series}
For any $T \in \calH_n^{nd}$, one has
\begin{eqnarray*}
b_\pi(T; t)  & = & \chi_\pi(\det T)^{n-1} \abs{\det T}^{t-n} 
\times 
\prod_{i=0}^{n-1}\, \frac{1 - (-1)^iq^{-t+i}}{1-(-1)^{i}q^{-(2n-t)+i}} \times 
b_\pi(T; 2n-t), 
\end{eqnarray*}
where $\chi_\pi(a) = (-1)^{v_\pi(a)}$ for $a \in k^\times$.
\end{thm}

\bigskip
\proof
Let us recall $\rho \in W$ given in Corollary~\ref{cor rho}. The value $s_t \in \C^n$ given by (\ref{s_t}) corresponds to $z^* \in \C^n$ in $z$-variable where $z_i^* = -\frac{t}{2}+i-\frac12 - (n-i+1)\frac{\pi\sqrt{-1}}{\log q}, \; 1 \leq i \leq n$, and $\rho(z^*)$ corresponds to 
$$
(1+\tfrac{\pi\sqrt{-1}}{\log q}, \ldots, 1\tfrac{\pi\sqrt{-1}}{\log q}) + (0, \ldots, 0, -\frac{t}{2}-\frac12 + (n-1)\tfrac{\pi\sqrt{-1}}{\log q})
$$
in $s$-variable.
By Corollary~\ref{cor rho}, we have
\begin{eqnarray*} \label{eq of omega-tilde}
\omega_T(x; s_t) =  \chi_\pi(\det T)^{n-1} \cdot \Gamma_\rho(z^*) \times \omega_T(x; s_{2n - t}).
\end{eqnarray*}
Hence we obtain by Proposition~\ref{prop2-1}, 
\begin{eqnarray} \label{junbi}
b_\pi(T; t) &=& \chi_\pi(\det T)^{n-1} 
\abs{\det T}^{t-n} \cdot  \gamma_n(t) \times b_\pi(T;2n-t).
\end{eqnarray}
where
\begin{eqnarray*}
\gamma_n(t) = \Gamma_\rho(z^*) \times \prod_{i=0}^{n-1} \frac{1 - q^{-t + 2i}}{1-q^{t-2(n-i)}} 
= \Gamma_\rho(z^*) \times (-1)^nq^{-nt+n(n+1)} \frac{1-q^{-t}}{1-q^{-t+2n}}.
\end{eqnarray*}
Since we have 
\begin{eqnarray*}
\Gamma_\rho(z^*) &=& 
\prod_{i<j}\, \frac{1 - (-1)^{i+j}q^{-t+i+j-2}}{(-1)^{i+j}q^{-t+i+j-1} - q^{-1}}
= 
(-q)^{\frac{n(n-1)}{2}} 
\prod_{i=1}^{n-1}\,\frac{1 - (-1)^{i}q^{-t+i}}{1 - (-1)^{n+i}q^{-t+n+i}}, 
\end{eqnarray*}
we get
\begin{eqnarray} 
\gamma_n(t)
&=&
(-1)^{\frac{n(n+1)}{2}}q^{-nt + \frac{n(3n+1)}{2}} 
\prod_{i=1}^{n}\, \frac{1-(-1)^{i-1} q^{-t+i-1}} {1-(-1)^{i+n}q^{-t+ n+i}} \nonumber \\
&=&
(-1)^{\frac{n(n+1)}{2}}q^{-nt + \frac{n(3n+1)}{2}} (-1)^{\frac{n(3n-1)}{2}}
q^{ns-\frac{n(3n+1)}{2}} 
\prod_{i=1}^{n}\, \frac{1-(-1)^{i-1} q^{-t+i-1}} {1-(-1)^{i+n}q^{-(2n-t)+ n- i}} \nonumber \\
&=&
\prod_{i=0}^{n-1}\, \frac{1-(-1)^{i} q^{-t+i}} {1-(-1)^{i}q^{-(2n-t)+i}}, 
\label{gam(t)}
\end{eqnarray}
hence we obtain the required functional equation of $b_\pi(T; t)$ by (\ref{junbi}).
\qed

\bigskip
\begin{rem} \label{rem: for decomp}
{\rm 
Let us recall the decomposition (\ref{shimura}) in Remark~\ref{shimura's decomp}.
Then by (\ref{gam(t)}), we see 
$$
\gamma_n(t) = f_T(q^{-t})/f_T(q^{t-2n}),
$$
and $\chi_\pi(\det T)^{n-1} \abs{\det T}^{t-n}$ gives the Gamma factor for the functional equation of $g_T(q^{-t})$.

The above functional equation is related to an element of the Weyl group of $U(H_n)$, which is not the case for (symmetric) Siegel series when $n$ is odd. F.~Sato and the author have studied in a similar line for Siegel series, we needed some harmonic analysis on $O(H_n)$ to establish the functional equations, and employed some previous results on particular $T$'s to determine the explicit Gamma factors. In the present case, we can obtain the explicit functional equations of hermitian Siegel series by a specialization of those of spherical functions $\omega_T(x; z)$. 
}
\end{rem} 

\bigskip
\begin{rem} \label{rem: Siegel series}
{\rm 
The existence of the functional equation of $b_\pi(T; t)$ was known in an abstract form as functional equations of Whittaker functions of $p$-adic groups by M.~L.~Karel \cite{Karel}. 
Recently T.~Ikeda \cite{Ikeda} has given explicit functional equations of $F_p(T; X) = g_T(X)$ on the basis of the results of S.~S.~Kudla and W.~J.~Sweet \cite{KS} for all quadratic extensions over $\Q_p$ containing split cases. There is a mistake in the range of $i$ of the definition of $t_p(K/\Q; X) = f_T(X)$ in \cite{Ikeda} p.1112, and it is better to refer the original $f_T(X)$ in \cite{Shimura} Theorem~13.6; if $K/\Q$ is unramified at $p$, $t_p(K/\Q; X)$ is the product of $1 - (-p)^i X$ from $i=0$ to $n-1$ as in Remark~\ref{shimura's decomp}.

}%
\end{rem}

\vspace{2cm}
\bibliographystyle{amsalpha}

\end{document}